\documentclass[american,
9pt,     
paper=a4,
oneside,
openany,
dvipsnames,
]{scrartcl}

\usepackage[
bottom=3.5cm
]{geometry}
\usepackage{amsmath,amssymb,amsfonts}%
\usepackage{amsthm}%
\usepackage{booktabs}
\usepackage{multirow,multicol}
\usepackage{caption}
\usepackage{diffcoeff}
\usepackage{mathtools}
\usepackage{hyperref}
\usepackage{cleveref}
\usepackage[ruled,longend,linesnumbered]{algorithm2e}
\usepackage{marvosym}
\usepackage{graphics}
\usepackage{graphicx}
\usepackage{blindtext}
\usepackage{wrapfig}
\usepackage{url}
\usepackage{tikz}
\usepackage{tkz-base,tkz-euclide}
\usepackage{tikz-3dplot}
\usetikzlibrary{calc,intersections,through,arrows.meta,external,shapes.geometric}
\usepackage{stmaryrd}
\usepackage{nicefrac}
\usepackage{listings}
\usepackage{xcolor}
\usepackage{geometry}
\usepackage{tabularx}
\usepackage{rotating}
\usepackage{subcaption}
\usepackage{enumitem}
\usepackage{lmodern}
\usepackage{placeins}

\usepackage[%
backend=biber, 
sorting=nyt, 
sortcites, 
safeinputenc, 
citestyle=numeric, %
bibstyle=ext-numeric, %
hyperref=true, 
maxbibnames=99, 
maxcitenames=2, 
url=false, 
isbn=false,
doi=true, 
giveninits=true,
date=year
]%
{biblatex}

\DeclareFieldFormat{titlecase:title}{\MakeSentenceCase*{#1}}

\renewbibmacro{in:}{%
  \ifentrytype{article}{}{\printtext{\bibstring{in}\intitlepunct}}}

\newtheoremstyle{bla}{}{}{\normalfont}{}{\bfseries}{}{0.5em}{}
\newtheoremstyle{bla2}{}{}{\itshape}{}{\bfseries}{}{0.5em}{}
\theoremstyle{bla2}
\newtheorem{theorem}{Theorem}[section]

\theoremstyle{bla}
\newtheorem{definition}[theorem]{Definition}%

\newcommand{\T}{\mathsf{T}}

\newcommand{\mop}{\Pi^{\mathsf{MO}}}

\newcommand{\lowi}{\mathcal{L}}

\newcommand{\oracle}{\mathsf{ALG}}

\begin{document}

\title{A Survey of Exact and Approximation Algorithms for Linear-Parametric Optimization Problems}

\author{Levin Nemesch\footnote{RPTU Kaiserslautern-Landau, Department of Mathematics, Paul-Ehrlich-Str. 14, 67663 Kaiserslautern, Germany, email:\{l.nemesch,stefan.ruzika\}@math.rptu.de}, Stefan Ruzika\footnotemark[\value{footnote}]~,    
     Clemens Thielen\footnote{Technical University of Munich, Campus Straubing for Biotechnology and Sustainability, Professorship of Optimization and Sustainable Decision Making, Am Essigberg~3, 94315~Straubing, Germany, email:\{clemens.thielen,alina.wittmann\}@tum.de}, Alina Wittmann\footnotemark[\value{footnote}]
}

\maketitle

\begin{abstract}
Linear-parametric optimization, where multiple objectives are combined into a single objective using linear combinations with parameters as coefficients, has numerous links to other fields in optimization and a wide range of application areas.
In this survey, we provide a comprehensive overview of structural results and algorithmic strategies for solving linear-parametric optimization problems exactly and approximately.
Transferring concepts from related areas such as multi-objective optimization provides further relevant results. The survey consists of two parts:
First, we list strategies that work in a general fashion and do not rely on specific problem structures.
Second, we look at well-studied parametric optimization problems and cover both important theoretical results and specialized algorithmic approaches for these problems.
Among these problems are parametric variants of shortest path problems, minimum cost flow and maximum flow problems, spanning tree problems, the knapsack problem, and matching problems.
Overall, we cover the results from 128 publications (and refer to 33 supplemental works) published between~1963 and~2024.

\end{abstract}


\section{Introduction}\label{subsec::Intro}

Linear-parametric optimization has been an important research area for more than 60~years with hundreds of articles published on the topic and a wide range of applications.

\medskip

In a linear-parametric optimization problem, different objectives are combined into a single objective by a linear combination with parameters as coefficients.
The goal is to find an optimal solution set, that is, a set containing an optimal solution for every possible combination of parameter values. 
This formulation naturally arises in various real-world applications where parameters change over time or across different scenarios.
In recent decades, the theory of parametric optimization has expanded significantly, along with its diverse applications in fields such as biomedical systems, power electronics, fleet planning, aeronautics, model predictive control, and process synthesis under uncertainty~\cite{Pistikopoulos_2012_Theoretical_and, oberdieck.diangelakis.ea2016MultiparametricProgrammingIts, Jenkins_1987_Using_Parametric}.
Furthermore, parametric optimization has many connections to other research areas, including computational geometry and various areas of optimization, such as multi-objective optimization, Lagrangian relaxation, robust optimization, sensitivity analysis, or multi-level optimization~\cite{helfrich.ruzika.ea2024EfficientlyConstructing,Fontaine_2014_Constraint-based_lagrangian,Avraamidou_2020_Adjustable_robust,Gal_1997_Advances_in,Avraamidou_2019_Multi-parametric_global,Tokuyama_2000_Combinatorics_on}.

By fixing the parameters of a parametric optimization problem, it becomes a conventional non-parametric optimization problem.
The parametric problem then is at least as hard as the resulting non-parametric problem.
But even if the non-parametric problem can be solved in polynomial time, the parametric problem itself might be intractable.
The reason for this is that, for many parametric problems, optimal solution sets can have a cardinality that is exponential in the encoding length of the problem input.
This includes, among others, the parametric variants of shortest path problems~\cite{Carstensen_1983_Complexity_of}, the assignment problem~\cite{Gassner_2010_A_fast}, and the minimum cost flow problem~\cite{Ruhe_1988_Complexity_results}.
Still, some problems such as the parametric minimum spanning tree problem remain tractable~\cite{Fernandez-Baca_1996_Using_sparsification}.

Both the facts that a parametric problem can be intractable and that it can be computationally hard for a fixed parameter value motivate the notion of approximation. 
Similarly to approximation in conventional optimization, approximation algorithms for parametric problems allow finding solution sets in polynomial time while ensuring a strict guarantee on the (worst-case) quality of the solutions~\cite{Holzhauser_2017_An_FPTAS, Giudici_2017_Approximation_schemes, Halman_2018_An_FPTAS,Bazgan_2022_An_approximation,Helfrich_2022_An_approximation}.

\medskip

This survey provides a comprehensive overview of theoretical developments, algorithmic strategies, and other recent advances in linear-parametric optimization, with a particular focus on some selected combinatorial optimization problems.
While existing overviews on parametric optimization are mostly concerned with parametric mixed integer linear programming~\cite{Geoffrion_1977_Exceptional_Paper,gal1995PostoptimalAnalysesParametric,greenberg1998AnnotatedBibliographyPostSolution,Pistikopoulos_2012_Theoretical_and,oberdieck.diangelakis.ea2016MultiparametricProgrammingIts}, our focus is on general, problem-agnostic approaches, approximation algorithms, and results for specific parametric problems.

\subsection{Scope and Outline}

This survey covers structural and algorithmic results for (multi-)parametric optimization problems.
In particular, cardinality results concerning optimal solution sets as well as exact solution methods and approximation algorithms for linear-parametric problems are addressed.

Parametric optimization is related to many other areas in optimization. Some of these relationships are briefly highlighted in Section~\ref{subsec::RelatedAreas}.
In Section~\ref{subsec::GeneralMethods}, we discuss general methods and results that apply to larger classes of parametric optimization problems.
We take advantage of the fact that some definitions and concepts in parametric optimization can also be found identically in multi-objective optimization.
Thus, we adopt some results from multi-objective optimization.
In Section~\ref{subsec::specificProblems}, we focus on methods and results for parametric variants of specific problems, in particular, well-known combinatorial optimization problems.
The selected problems are among the most-studied problems in the context of parametric optimization throughout peer-reviewed journals and conference proceedings in the English language.

\subsection{Preliminaries and Definitions}

We distinguish between the natural numbers~$\mathbb{N}$ starting with~$1$ and~$\mathbb{N}_0$ starting with~$0$.
Given $p$, $n \in \mathbb{N}$, let $[p]\coloneq\{1,\dots,p\}$ be the set of natural numbers up to~$p \in \mathbb{N}$.
We denote the nonnegative orthant by $\mathbb{R}^p_{\geq 0} \coloneq \{y \in \mathbb{R}^p: \, y \geq 0 \}$, where $0 \in \mathbb{R}^p$ is the $p$-dimensional zero vector and $\geq$ is the weak component-wise order defined by
\begin{align*}
    y \geq y' \text{ if and only if } y_i \geq y_i' \quad \forall i \in [p].
\end{align*}
Given a function $f:\mathbb{R}^p\rightarrow \mathbb{R}^n$ and a set $Z \subseteq \mathbb{R}^p$, we denote the image of~$Z$ under~$f$ by $f(Z)\coloneq \{f(z): \, z \in Z\} \subseteq \mathbb{R}^n$.

\medskip

For running time specifications, we use the well-known Landau notation.
We also use $\tilde{\mathcal{O}}(...)$ (read soft-O), an extension to the Landau notation that hides polylogarithmic factors, e.g., $\tilde{\mathcal{O}}(g(n))$ for $\mathcal{O}(g(n)\log^k(n))$ for some $k \in \mathbb{N}$.
For a problem~$\mathsf{P}$, $T_{\mathsf{P}}$ denotes the running time of an oracle that solves~$\mathsf{P}$.
A polynomial in~$n$ is referred to as $\operatorname{poly}(n)$.
An algorithm is \emph{(input-)polynomial} if, for an input of size~$I$, it has a running time in~$\mathcal{O}(\operatorname{poly}(I))$.
If no polynomial time algorithm exists for a problem, it is called  \emph{intractable}.
An algorithm is \emph{output-polynomial} if, for an input of size~$I$ and an output of size~$O$, it has a running time in~$\mathcal{O}(\operatorname{poly}(I,O))$.
\medskip

The following definitions introduce general (multi-)parametric minimization problems~\cite{Bazgan_2022_An_approximation,Helfrich_2022_An_approximation}.
Note that, for the sake of simplicity, we consider minimization problems, but all concepts can be defined analogously for maximization problems.

\begin{definition}[Parameter Set]\label{def:ParamSet}
A \emph{parameter set~$\Lambda\subseteq\mathbb{R}^p$} is a set of the form
    \[
    \Lambda \coloneq \left\{\lambda\in\mathbb{R}^p: \lambda_{\min}\leq\lambda\leq\lambda_{\max} \right\}
    \]
    for a given $\lambda_{\min}\in\left(\mathbb{R}\cup \{-\infty\}\right)^p$ and a given $\lambda_{\max}\in\left(\mathbb{R}\cup \{\infty\}\right)^p$.
\end{definition}

We remark that parametric optimization problems can easily be defined for more complex parameter sets, e.g., for polyhedral or non-convex parameter sets.
Almost all parameter sets used in the literature, however, are covered by \Cref{def:ParamSet}.

\begin{definition}[(Multi-)Parametric Minimization Problem]
Let $p\in \mathbb{N}$.
A \emph{linear \mbox{($p$-)parametric} minimization problem} $\Pi=(X, f, \Lambda)$ is a triple consisting of
\begin{enumerate}[label=(\alph*)]
    \item a nonempty set~$X\subseteq\mathbb{R}^n$ of feasible solutions for some $n\in\mathbb{N}$,
    \item a function  $f\colon X \to \mathbb{R}^{p+1}, \, x \mapsto f(x) = (f_0(x),f_1(x),\dots,f_p(x))^\mathsf{T}$, and
    \item a \emph{parameter set} $\Lambda\in\mathbb{R}^p$.
\end{enumerate}
A solution of~$\Pi$ is a set~$S\subseteq X$ that contains, for every $\lambda\in\Lambda$, an optimal solution of the problem
    \begin{align*}\tag{$\Pi(\lambda)$}
        \min_{x\in X}F_\lambda(x)\coloneq\min_{x\in X} f_0(x) + \sum_{i=1}^p \lambda_i f_i(x).
    \end{align*}
\end{definition}

For a fixed parameter value~$\lambda\in\Lambda$, the problem~$\Pi(\lambda)$ is a conventional (non-parametric) minimization problem.
The \emph{optimal cost curve} $F^*:\Lambda \to \mathbb{R}\cup\{-\infty\}$, $\lambda \mapsto F^*(\lambda)\coloneq\inf_{x \in X} F_{\lambda}(x)$ assigns to every $\lambda\in\Lambda$ the optimal objective value of~$\Pi(\lambda)$, or $-\infty$ if $\Pi(\lambda)$ is unbounded.
The set~$F^*(\Lambda)$ is called the \emph{lower envelope} of~$\Pi$.

With $\oracle$, we denote an oracle for a problem~$\Pi$ that takes as input any $\lambda\in\Lambda$ and computes an optimal solution of $\Pi(\lambda)$.

For an optimal solution set~$S$, we do not only expect it to contain an optimal solution for every~$\lambda\in\Lambda$, but to also have minimal cardinality among all such sets.
With some slight abuse of notation, $B$ denotes the cardinality $\left|S \right|$ if it is finite or $\infty$ if it is infinite.

The \emph{critical region}~$\Lambda(x)\subseteq\Lambda$ of a solution~$x\in X$ is the set of all parameters for which $x$ is an optimal solution of $\Pi(\lambda)$.
In the case of a single-parametric problem, the critical region of each solution in an optimal solution set is an interval in~$\mathbb{R}$.
The point where two adjacent intervals intersect is called a \emph{breakpoint}.

As the cardinality of an optimal solution set can be superpolynomially large in the worst case, even if $p = 1$, approximation provides a good concept to substantially reduce the number of required solutions while still obtaining provable
solution quality.

\begin{definition}
    For an approximation factor~$\alpha \geq 1$ and a fixed parameter vector~$\lambda \in \Lambda$, a feasible solution $x \in X$ is called \emph{$\alpha$-approximate} for the problem~$\Pi(\lambda)$ if $F_{\lambda}(x) \leq \alpha \cdot F_{\lambda}(x')$ for all $x' \in X$.
\end{definition}

A necessary condition for the existence of an $\alpha$-approximate solution of~$\Pi(\lambda)$ is that the optimal objective value of~$\Pi(\lambda)$ is greater than or equal to zero.
Therefore, it is assumed that~$F^*(\lambda)\geq0$ for all~$\lambda\in\Lambda$ whenever parametric approximation algorithms are considered.

Usually, the goal in parametric approximation is to find a set of solutions that contains an~$\alpha$-approximate solution for each non-parametric problem~$\Pi(\lambda)$:

\begin{definition}[Approximation set] \label{def:ApproxSet}
    For an approximation factor $\alpha \geq 1$, a finite set $S \subseteq X$ is called an \emph{$\alpha$-approximation set} to~$\Pi$ if it contains an $\alpha$-approximate solution $x \in S$ of $\Pi(\lambda)$ for every $\lambda \in \Lambda$.
    An $\alpha$-approximation set is \emph{minimal} if removing any solution would no longer provide an $\alpha$-approximation set.

    An algorithm $A$ that computes an $\alpha$-approximation set for an instance~$\Pi$ in time polynomially bounded in the instance size is called an \emph{$\alpha$-approximation algorithm}.

    A \emph{polynomial-time approximation scheme} (PTAS) is a family of algorithms $(A_{\varepsilon})_\varepsilon$ such that, for any $\varepsilon > 0$, the algorithm $A_{\varepsilon}$ is a $(1+\varepsilon)$-approximation algorithm.
    A PTAS~$(A_{\varepsilon})_\varepsilon$ is a \emph{fully polynomial-time approximation scheme} (FPTAS) if the running time of~$A_{\varepsilon}$ is also polynomial in $\nicefrac{1}{\varepsilon}$.
\end{definition}

Similar to the oracle $\oracle$, an \emph{approximation oracle} $\oracle_\beta$ for a $\beta\geq1$ takes as input any $\lambda\in\Lambda$ with $F^*(\lambda)\geq0$ and returns a $\beta$-approximate solution of the problem $\Pi(\lambda)$.


\section{Related Areas}\label{subsec::RelatedAreas}

We relate parametric optimization to some other selected areas of mathematical optimization. The purpose is twofold: We demonstrate the impact of parametric optimization and, vice versa, we illustrate how parametric optimization may benefit from progress made in related fields of research.

\medskip

In Lagrangian relaxation, constraints are moved into the objective function together with a Lagrangian multiplier such that any violation of the constraints is penalized~\cite{Marshall_1981_The_Lagrangian}.
\begin{align*}
    \begin{array}{l r l}
        {\displaystyle \min_{x}}\quad{} & f_0(x) \\
        \text{s.t.} & f_1(x) & \leq 0 \\
        & g(x) & \leq 0 \\
    \end{array}
    \qquad \rightarrow \qquad
    \left\{ \begin{array}{l l}
        \displaystyle\min_{x}\quad{} & f_0(x) + \lambda\cdot f_1(x) \\
        \text{s.t.}  & g(x) \leq 0 \\
    \end{array} \right\}_{\lambda \in \Lambda}
\end{align*}
The degree of penalization depends on the Lagrangian multiplier which influences the quality of the dual bound obtained.
Lagrangian relaxation and solving the Lagrangian dual problem can be viewed as a parametric optimization problem (see~\cite{Bryson_1991_Parametric_programming,Bryson_1993_A_parametric,Fontaine_2014_Constraint-based_lagrangian}).

\medskip

In many domains such as engineering, computer science, or finance, data of related to the objective function and / or the constraints of optimization problems is often uncertain. Robust optimization aims to find solutions that are robust with respect to uncertainty described by realizations of parameter values within a given set of scenarios. 
Solution approaches for robust optimization problems are developed with the help of parametric optimization by considering the uncertainty set as the parameter set, which means that solving the parametric problem provides a solution for every possible scenario in the uncertainty set~\cite{Gabrel_2014_Recent_advances,Chuong_2019_Exact_relaxations,Avraamidou_2020_Adjustable_robust,Antczak_2021_Parametric_Approach}.

\medskip

In the field of sensitivity analysis, the behavior of primal solutions is examined when there is perturbation or uncertainty in the input parameters or the structure of an optimization problem.
This means that results and methods for parametric optimization problems can also be found in the literature that is concerned with sensitivity analysis for combinatorial optimization problems (see~\cite{Gusfield_1980_Sensitivity_Analysis,Gal_1986_Stability_in}).
An overview on the historical development and literature overview of sensitivity analysis and parametric optimization can be found in~\cite{Gal_1997_Advances_in}.

\medskip

A multi-level optimization problem consists of several optimization problems embedded in each other, where the constraints of the problem on one level are implicitly determined by the problems on the lower levels.
Multi-parametric optimization is used as a tool for solving multi-level optimization problems (see~\cite{Dempe_2020_Bilevel_Optimization,Pistikopoulos_2012_Theoretical_and,Avraamidou_2019_Multi-parametric_global,Goyal_2021_Parametric_approach,Avraamidou_2022_Multi-level_Mixed}).
The idea is to consider the lowest-level problem as a multi-parametric problem with the higher-level decision variables as parameters. 
Then, this parametric problem is solved, and each critical region defines a new multi-level problem having one level less.
Iteratively, this process is continued until one single level remains~\cite{Avraamidou_2019_Multi-parametric_global}.

\medskip

\citeauthor{Tokuyama_2000_Combinatorics_on}~\cite{Tokuyama_2000_Combinatorics_on} compares the trajectories of all weight functions of a parametric optimization problem in a plane with line arrangements in computational geometry. Apart from line arrangements, lower envelopes of line segments and $k$-sets are also connected to parametric optimization~\cite{Eppstein_1998_Geometric_Lower}.

\subsection{Multi-Objective Optimization}\label{subsec::moopt}

An area with particularly close connections to parametric optimization is \emph{multi-objective optimization}.
We only outline some key concepts here. For an in-depth introduction to multi-objective optimization, we refer to \citeauthor{ehrgott2005multicriteriaOptimization}~\cite{ehrgott2005multicriteriaOptimization}.

\medskip

Given a set~$X\in\mathbb{R}^n$ of feasible solutions and a (vector-valued) objective function~$f=(f_1,\dots,f_d)\colon X\rightarrow \mathbb{R}^d$ for $d\in\mathbb{N}_{>1}$, a \emph{$d$-objective optimization problem} is of the form:
\begin{align*}\tag{$\mop$}
    \min_{x\in X} & \left(f_1(x),\ldots , f_{d}(x)\right)
\end{align*}
Since solutions optimizing $f_1,\dots,f_d$ simultaneously do typically not exist for a $d$-objective problem, one is usually interested in so-called \emph{efficient solutions}~\cite{ehrgott2005multicriteriaOptimization}. Here, a feasible solution~$x\in X$ is called \emph{efficient} (or \emph{Pareto optimal}) if there is no other feasible solution~$x'\in X$ such that $f(x')\leq f(x)$ and $f(x') \neq f(x)$. The image~$f(x)$ of an efficient solution~$x$ under~$f$ is called a \emph{non-dominated} image. The goal in a multi-objective optimization problem then typically consists of determining the set of all non-dominated images (the \emph{non-dominated set}) and one corresponding efficient solution for each of these non-dominated images.

\medskip

A significant part of the multi-objective optimization literature, however, is mainly interested in a subset of the non-dominated set, namely in the set of non-dominated extreme points of the convex hull $\operatorname{conv}(Y)$ of the image set~$Y=f(X)$. 
This is where the close relation between multi-objective and parametric optimization arises.
The set of non-dominated extreme points of~$f(X)$ can also be characterized as an optimal solution set of the parametric problem $\min_{x\in X} \lambda^\T f(x)$ with the parameter set $\Lambda = \mathbb{R}_{\geq0}^d$~\cite{ehrgott2005multicriteriaOptimization}.
In the multi-objective literature, this parameter set is called the \emph{weight set}.
Hence, solving this parametric problem and finding the non-dominated extreme points are equivalent problems.

\medskip

This relation led to the development of some similar approaches in the parametric and the multi-objective optimization literature.
A classic example for this is the famous technique known as the \emph{Eisner-Severance method} in parametric optimization and as the \emph{dichotomic approach} in multi-objective optimization (see \Cref{subsec::ESmethod}).
In \Cref{subsec::decomplowi}, we provide an overview of some algorithmic strategies that were independently developed in both areas and whose connection has, to the best of our knowledge, not previously been outlined in the literature.

\medskip

The similarities between parametric and multi-objective optimization also allow approximation algorithms from the multi-objective setting to be used in the parametric setting.
As long as $\Lambda=\mathbb{R}^d_{\geq0}$, every so-called \emph{$\alpha$-convex Pareto set} for the corresponding multi-objective problem is an $\alpha$-approximation set of the parametric problem, and vice versa~\cite{helfrich.ruzika.ea2024EfficientlyConstructing}.
Under the same condition, every so-called \emph{$\alpha$-approximate Pareto set} (a stricter notion of approximation in the multi-objective setting) is also a parametric $\alpha$-approximation set, but it might contain solutions that are redundant for the latter.
These can be filtered out by convex hull algorithms~\cite{Giudici_2017_Approximation_schemes}.
For a comprehensive overview on multi-objective approximation, we refer to the recent survey in~\cite{herzel.ruzika.ea2021ApproximationMethodsMultiobjective}.


\section{General Methods}\label{subsec::GeneralMethods}

In this section, we present general exact and approximate methods that are applicable for a broad class of parametric optimization problems, independent of specific problem structures. The Eisner-Severance method is discussed in Section \ref{subsec::ESmethod}, Gusfield’s algorithm using Megiddo’s method in Section \ref{subsec::megiddogusfield}, the connection between parametric and multi-objective optimization in Section \ref{subsec::decomplowi}, parametric approximation in Section \ref{subsec::generalAPX} and results for parametric mixed integer linear programs are reviewed in Section \ref{subsec::milp}.

\subsection{The Eisner-Severance Method}\label{subsec::ESmethod}

A widely-used algorithm is the \emph{Eisner-Severance method~\cite{Eisner_1976_Mathematical_Techniques}}, which has also been introduced by the multi-objective literature under a different name.
\citeauthor{Eisner_1976_Mathematical_Techniques} aimed at record-segmentation problems in shared databases, but their method can be used as a general algorithm to solve single-parametric optimization problems.
The precise origin of the method is unclear, as \citeauthor{Nauss_1975_Parametric_Integer}~\cite{Nauss_1975_Parametric_Integer} independently developed a similar approach for Knapsack problems in 1975.
\citeauthor{Radke_1975_Sensitivity_Analysis}~\cite{Radke_1975_Sensitivity_Analysis} and \citeauthor{Geoffrion_1977_Exceptional_Paper}~\cite{Geoffrion_1977_Exceptional_Paper} both generalize \citeauthor{Nauss_1975_Parametric_Integer}' algorithm to more generic problems and provide proofs of correctness and running time.
\citeauthor{Gusfield_1980_Sensitivity_Analysis}~\cite{Gusfield_1980_Sensitivity_Analysis} combines both branches of the literature and states a refined proof of correctness and running time.

\medskip

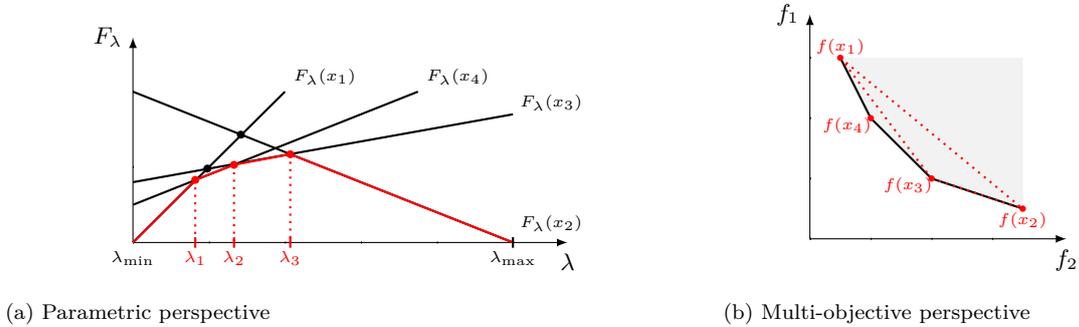
\begin{figure}[h!]
    \centering
    \begin{minipage}[b]{0.6\textwidth} 
	\centering
	\begin{tikzpicture}[scale=1]
		\tkzInit[xmax=5.22,ymax=2.22,xmin=0,ymin=0]
		\tkzSetUpAxis[ticka=0pt, tickb=0pt]
		\tkzDrawXY[/tkzdrawX/label=$\lambda$,/tkzdrawY/label=$F_\lambda$]
		
		\draw [thick] (0,0)--(2,2);
		\draw [thick] (0,2)--(5,0);
		\draw [thick] (0,0.5)--(3.75,2);
		\draw [thick] (0,0.8)--(5,1.7);

		\draw [red,thick] (0,0)--(0.82,0.83);
		\draw [red,thick] (0.82,0.83)--(1.33,1.03);
		\draw [red,thick] (1.33,1.03)--(2.07,1.17);
		\draw [red,thick] (2.07,1.17)--(5,0);
		
		\fill [red](0.82,0.83) circle(1.5pt);
		\fill [red](1.33,1.03) circle(1.5pt);
		\fill [red](2.07,1.17) circle(1.5pt);
		\fill [black](1.42,1.43) circle(1.5pt);
		\fill [black](0.98,0.98) circle(1.5pt);
		
		\draw [red,dotted,thick] (0.82,0.83)--(0.82,0);
		\draw [red,dotted,thick] (1.33,1.03)--(1.33,0);
		\draw [red,dotted,thick] (2.07,1.17)--(2.07,0);
		
		\draw [red,thick] (0.82,0.07)--(0.82,-0.07);
		\draw [red,thick] (1.33,0.07)--(1.33,-0.07);
		\draw [red,thick] (2.07,0.07)--(2.07,-0.07);
		\draw [thick] (5,0.07)--(5,-0.07);
		
		\tkzText[below,red](0.82,0){\scriptsize $\lambda_1$}
		\tkzText[below,red](1.33,0){\scriptsize $\lambda_2$}
		\tkzText[below,red](2.07,0){\scriptsize $\lambda_3$}
		\tkzText[below](5,0){\scriptsize $\lambda_{\text{max}}$}
		\tkzText[below](0,0){\scriptsize $\lambda_{\text{min}}$}
		
		\tkzText[right](2,2.2){\scriptsize $F_{\lambda}(x_1)$}
		\tkzText[right](4.98,0.25){\scriptsize $F_{\lambda}(x_2)$}
		\tkzText[right](4.98,1.85){\scriptsize $F_{\lambda}(x_3)$}
		\tkzText[right](3.75,2.2){\scriptsize $F_{\lambda}(x_4)$}
	\end{tikzpicture}
\end{minipage}
\hfill
\begin{minipage}[b]{.35\textwidth} 
	\centering
	\begin{tikzpicture}[scale=0.8]
		\tkzInit[xmax=3.72,ymax=3.22,xmin=0,ymin=0]
		\tkzSetUpAxis[ticka=0pt, tickb=0pt]
		\tkzDrawXY[/tkzdrawX/label=$f_2$,/tkzdrawY/label=$f_1$]
		
		\draw [thick] (1,2)--(2,1);
		\draw [thick] (1,2)--(0.5,3);
		\draw [thick] (2,1)--(3.5,0.5);
		
		\fill [red](3.5,0.5) circle(1.5pt);
		\fill [red](0.5,3) circle(1.5pt);
		\fill [red](1,2) circle(1.5pt);
		\fill [red](2,1) circle(1.5pt);
		
		\begin{scope}[on background layer]
			\fill[gray!10] (0.5,3)--(1,2)--(2,1)--(3.5,0.5)--(3.5,3);
		\end{scope}
		
		\draw [red,dotted,thick] (0.5,3)--(3.5,0.5);
		\draw [red,dotted,thick] (2,1)--(0.5,3);
		\draw [red,dotted,thick] (2,1)--(3.5,0.5);
		
		\tkzText[above=-2pt,red](0.5,3){\scriptsize $f(x_1)$}
		\tkzText[below left=-5pt,red](1,2){\scriptsize $f(x_4)$}
		\tkzText[below left=-5pt,red](2,1){\scriptsize $f(x_3)$}
		\tkzText[below=-2pt,red](3.5,0.5){\scriptsize $f(x_2)$}
	\end{tikzpicture}
\end{minipage}
\begin{minipage}[t]{0.6\textwidth}
		\centering
    \subcaption{Parametric perspective}\label{fig::ES1}
  \end{minipage}
  \hfill
  \begin{minipage}[t]{0.35\textwidth}
  	\centering
    \subcaption{Multi-objective perspective}\label{fig::ES2}
  \end{minipage}
    \caption{Example showing two different perspectives on the Eisner-Severance method~\cite{Eisner_1976_Mathematical_Techniques}, returning the optimal cost curve in the parametric perspective in (a) and the non-dominated extreme points in the multi-objective perspective in (b).}
    \label{Fig:ESmethod}
\end{figure}

The Eisner-Severance method works on single-parametric problems with parameter sets of the form $[\lambda_{\min},\lambda_{\max}]$.
It requires access to an oracle $\oracle$.
The method is illustrated in Figure~\ref{Fig:ESmethod}(a) from the parametric perspective and in Figure~\ref{Fig:ESmethod}(b) from the multi-objective perspective.
After initializing $\{\lambda_{\text{min}},\lambda_{\text{max}}\}$, two solutions $x_1\coloneq \oracle(\lambda_{\text{min}})$ and $x_2\coloneq \oracle(\lambda_{\text{max}})$ are computed.
Then, the parameter value~$\lambda^*$ is computed for which the objective function value of both solutions is equal, i.e., the point where the functions $F_{\lambda^*}(x_1)$ and $F_{\lambda^*}(x_2)$ intersect in Figure~\ref{Fig:ESmethod}(a) or the normal of the connecting line between~$f(x_1)$ and~$f(x_2)$ in Figure~\ref{Fig:ESmethod}(b).
If $F_{\lambda^*}(x_1)$ is the optimal objective value at~$\lambda^*$, then~$\lambda^*$ is a breakpoint, $x_1$ is optimal in $[\lambda_{\min},\lambda^*]$ and $x_2$ is optimal in $[\lambda^*,\lambda_{\max}]$. Thus, $x_1$ and $x_2$ together are an optimal solution set in the interval $[\lambda_{\min},\lambda_{\max}]$.
Otherwise, the above procedure is applied recursively to the two sub-intervals $[\lambda_{\min},\lambda^*]$ and $[\lambda^*,\lambda_{\max}]$.

\medskip

Since all other operations can be performed in~$\mathcal{O}(1)$, the number of calls to~$\oracle$ determines the running time of the Eisner-Severance method.
It is called at least $2B-1$ and at most $3B-2$ times (cf.~\cite{csirmaz2021InnerApproximationAlgorithm}).
Under the additional assumption that~$\oracle$ only outputs solutions that correspond to line segments of the optimal cost curve, exactly $2B-1$ calls to~$\oracle$ are needed~\cite{Aneja_1979_Bicriteria_Transportation, Gusfield_1980_Sensitivity_Analysis}.

\medskip

An approximation algorithm for single-parametric problems based on the Eisner-Severance method is developed by \citeauthor{Bazgan_2022_An_approximation}~\cite{Bazgan_2022_An_approximation}.
If a $\beta$-approximation oracle $\oracle_\beta$ for $\Pi(\lambda)$ and all $\lambda\in\Lambda$ is available, their algorithm is able to compute a $(1+\varepsilon)\beta$-approximation set for any~$\varepsilon>0$.
Both the cardinality of the approximation set and the number of calls to~$\oracle_\beta$ are bounded by a polynomial in the size of the input and~$\nicefrac{1}{\varepsilon}$.
\Cref{Fig:ESmethod_Bazgan} illustrates this approach.

\begin{figure}[h!]
    \begin{minipage}{\linewidth}
        \centering
        \begin{tikzpicture}[scale=1.2]
    \tkzInit[xmax=5.22,ymax=2.22,xmin=0,ymin=0]
    \tkzSetUpAxis[ticka=0pt, tickb=0pt]
    \tkzDrawXY[/tkzdrawX/label=$\lambda$,/tkzdrawY/label=$F_{\lambda}$]

    \draw [red,dotted,thick] (0.82,0.83)--(0.82,0);
    \draw [dotted,thick] (1.42,0.92)--(1.42,0);
    \draw [red,dotted,thick] (3.77,1.67)--(3.77,0);
    
    \draw [thick] (0,0)--(2,2);
    \draw [thick] (0,1.28)--(5,1.8);
    \draw [red, thick] (0,0.59)--(4.9,2);
    \draw [gray, thick] (0,0.89)--(3.8,2); 
    \draw [thick] (0,0.7)--(4.3,2);

    
    \fill [red](0.82,0.83) circle(1.5pt);
    \fill [gray](1.42,1.298) circle(1.5pt);
    \fill [red](1.42,1) circle(1.5pt);
    \fill [red](3.77,1.67) circle(1.5pt);
    \fill [black](1.42,1.43) circle(1.5pt);
    
    \draw [red,thick] (0.82,0.07)--(0.82,-0.07);
    \draw [thick] (1.42,0.07)--(1.42,-0.07);
    \draw [red,thick] (3.77,0.07)--(3.77,-0.07);
    \draw [thick] (5,0.07)--(5,-0.07);
    
    \tkzText[below,red](0.82,0){\scriptsize $\lambda_{\text{l}}$}
    \tkzText[above,red](5,2){\scriptsize $R(\lambda)$}
    \tkzText[above,gray](3.7,2){\scriptsize $(1+\varepsilon)R(\lambda)$}
    \tkzText[below](1.42,0){\scriptsize $\lambda_{\text{m}}$}
    \tkzText[below,red](3.77,0){\scriptsize $\lambda_{\text{r}}$}
    \tkzText[below](5,0){\scriptsize $\lambda_{\text{max}}$}
    \tkzText[below](0,0){\scriptsize $\lambda_{\text{min}}$}
    
\end{tikzpicture}
        \caption*{\footnotesize The algorithm starts with a parameter interval $([\lambda_{\text{l}},\lambda_{\text{r}}],\oracle_\beta(\lambda_{\text{l}}),\oracle_\beta(\lambda_{\text{r}}))$. First, it is checked whether the optimal / approximate solution~$x_{\text{l}}\coloneq\oracle_\beta(\lambda_{\text{l}})$ for $\lambda_{\text{l}}$ or~$x_{\text{r}}\coloneq\oracle_\beta(\lambda_{\text{r}})$ for $\lambda_{\text{r}}$ is $\alpha$-approximate for the whole interval $[\lambda_{\text{l}},\lambda_{\text{r}}]$ by examining the two conditions $F_{\lambda_{\text{r}}}(x_{\text{l}}) \leq F_{\lambda_{\text{r}}}(x_{\text{r}})$ and $F_{\lambda_{\text{l}}}(x_{\text{r}}) \leq F_{\lambda_{\text{l}}}(x_{\text{l}})$. 
        Otherwise, it is checked whether the intersection point $(\lambda_{\text{m}},y_{\text{m}})$ of $F_\lambda(x_{\text{l}})$ and $F_\lambda(x_{\text{r}})$ lies below the red connecting line~$R(\lambda)$ shifted by a factor of~$1+\varepsilon$, i.e., below the gray line~$(1+\varepsilon)R(\lambda)$. Otherwise, the algorithm is recursively called with $\left([\lambda_{\text{l}},\lambda_{\text{m}}],x_{\text{l}},\oracle_\beta(\lambda_{\text{m}})\right)$ and $\left([\lambda_{\text{m}},\lambda_{\text{r}}],\oracle_\beta(\lambda_{\text{m}}),x_{\text{r}}\right)$.}
    \end{minipage}
    \caption{Approximate variant of the Eisner-Severance method from \cite{Bazgan_2022_An_approximation} \label{Fig:ESmethod_Bazgan}}
\end{figure}

A strategy that works functionally identical to the Eisner-Severance method was independently developed in multi-objective optimization to identify non-dominated extreme points of bi-objective optimization problems.
There, it is known as the \emph{dichotomic approach} and first appears in~\cite{Cohon_1978_Multiobjective_programming,Aneja_1979_Bicriteria_Transportation,Dial_1979_A_model}.
An approximate variant of this algorithm also exists. It is known as the \emph{Sandwich Algorithm with Underlying Chord Partition Rule} \cite{Burkard_1991_Sandwich_approximation,Rote_1992_The_convergence} or the \emph{Chord Algorithm} \cite{Diakonikolas_2011_Approximation_of,Daskalakis_2016_How_Good}, and is illustrated in \Cref{Fig:ChordMethod}.

\begin{figure}[h!]
    \begin{minipage}{\linewidth}
        \centering
        \begin{tikzpicture}[scale=1]
	\tkzInit[xmax=3.72,ymax=3.22,xmin=0,ymin=0]
	\tkzSetUpAxis[ticka=0pt, tickb=0pt]
	\tkzDrawXY[/tkzdrawX/label=$f_2$,/tkzdrawY/label=$f_1$]
	
	\fill [gray](0.8,0.8) circle(1.5pt);
	
	\begin{scope}[on background layer]
		\fill[gray!10] (0.5,3)--(1,2)--(2,1)--(3.5,0.5)--(3.5,3);
	\end{scope}
	
	\draw [red,dotted,thick] (0.5,3)--(3.5,0.5);
	\draw [red,dotted,thick] (2,1)--(0.5,3);
	\draw [red,dotted,thick] (2,1)--(3.5,0.5);
	
	\draw [gray,semithick] (0.8,0.8)--(3.5,0.5);
	\draw [gray,semithick] (0.8,0.8)--(0.5,3);
	\draw [red,semithick] (2.6,0.5)--(0.5,2.25);
	
	\fill [red](2.47,0.61) circle(1.5pt);
	\fill [red](0.615,2.15) circle(1.5pt);
	
	\draw [thick] (1,2)--(2,1);
	\draw [thick] (1,2)--(0.5,3);
	\draw [thick] (2,1)--(3.5,0.5);
	
	\fill [black](3.5,0.5) circle(1.5pt);
	\fill [black](0.5,3) circle(1.5pt);
	\fill [black](1,2) circle(1.5pt);
	\fill [black](2,1) circle(1.5pt);
	
	\tkzText[above,black](0.5,3){\scriptsize $f(x_1)$}
	\tkzText[below left=-3pt,gray](0.8,0.8){\scriptsize $s_1$}
	\tkzText[left=-1pt,red](0.5,2.25){\scriptsize $s_2$}
	\tkzText[below=0pt,red](2.6,0.5){\scriptsize $s_3$}
	\tkzText[left,black](2,0.95){\scriptsize $f(x_3)$}
	\tkzText[below,black](3.5,0.5){\scriptsize $f(x_2)$}
\end{tikzpicture}
        \caption*{\footnotesize Again, the aim is to find a solution that is optimal with respect to the normal of the connecting line between $f(x_1)$ and $f(x_2)$.
        Let the two gray lines be the connecting lines from the previous iterations shifted to the optimal points.
        In contrast to the exact Eisner-Severance method, the recursion already terminates if the connecting line between $f(x_1)$ and $f(x_2)$ is not more than the approximation factor away from the intersection point $s_1$ or the new optimal / approximate image $f(x_3)$.}
    \end{minipage}
    \caption{Chord Algorithm}
    \label{Fig:ChordMethod}
\end{figure}

\subsection{Gusfield's Algorithm and Megiddo's Method}\label{subsec::megiddogusfield}

Another general algorithm for single-parametric optimization problems is described by \citeauthor{Gusfield_1983_Parametric_Combinatorial}~\cite{Gusfield_1983_Parametric_Combinatorial}.
It is based on a well-known technique from single-parametric search, called \emph{Megiddo's method}.
First described by \citeauthor{megiddo1979CombinatorialOptimizationRational}~\cite{megiddo1979CombinatorialOptimizationRational} for parametric search in combinatorial optimization problems, Megiddo's method is developed further in~\cite{megiddo1983ApplyingParallelComputation,cole1987SlowingSortingNetworks,vanoostrum.veltkamp2004ParametricSearchMade,cohen.megiddo1993MaximizingConcaveFunctions}, adapted for non-linear problems in~\cite{toledo1993MaximizingNonLinearConcave,agarwal.sharir.ea1994ApplicationsParametricSearching}, and turned into an approximation algorithm in~\cite{hashizume.fukushima.ea1987ApproximationAlgorithmsCombinatoriala,toledo1993ApproximateParametricSearching}.

\medskip

The goal of Megiddo's method is to find a parameter value~$\lambda^*$ for which some condition holds -- usually a parameter value where~$F^*$ attains its global optimum.
As long as a polynomial-time algorithm for~$\Pi(\lambda)$ is available for every $\lambda\in\Lambda$, Megiddo's method also runs in polynomial time.
The basic idea is to simulate the algorithm for~$\Pi(\lambda^*)$ without actually knowing the value of~$\lambda^*$.
Whenever the simulated algorithm encounters a conditional branching, Megiddo's method needs to detect which branch the algorithm would take for the (still unknown) input~$\lambda^*$.

\medskip

Any conditional branching that is encountered is a less-or-equal comparison between two linear functions of~$\lambda$ (or can be decided trivially).
The result of the comparison is determined by the parameter value of~$\lambda^*$.
However, it can also be decided by solving~$\Pi(\lambda)$ for the parameter value for which the two functions intersect.
Solving~$\Pi(\lambda)$ can be done with the same algorithm that is simulated for~$\lambda^*$.
When the simulated algorithm terminates, the value of~$\lambda^*$ and a solution of~$\Pi(\lambda^*)$ are known.

\medskip

Gusfield's algorithm transforms this parametric search technique into an algorithm that computes an optimal solution set.
It uses Megiddo's method to find an initial breakpoint, and then to iteratively identify subsequent breakpoints with corresponding solutions.
The running time of Gusfield's algorithm is in $\mathcal{O}(B\cdot T_{\Pi}^2)$~\cite{Gusfield_1983_Parametric_Combinatorial}.

\subsection{Minimization Diagram and Lower Image}\label{subsec::decomplowi}

Next, we look at some general strategies that are applicable to every multi-parametric mixed integer linear program that is bounded in $\Lambda$.
We also include some algorithms from multi-objective optimization that, fundamentally, take the parametric perspective to find non-dominated extreme points.
To keep the notation simple, we directly transfer all results from the multi-objective literature to the parametric setting.

\subsubsection{Decomposition of the Parameter Set and Minimization Diagrams}\label{subsec::decomposition}

For every~$x\in X$ where $X$ is described by a system of mixed integer linear constraints, the critical region~$\Lambda(x)$ is a convex polyhedron.
Several approaches decompose the parameter set~$\Lambda$ into disjoint $p$-dimensional critical regions and find one solution per critical region.
Such a subdivision is called a \emph{minimization diagram}~\cite{schwartz.sharir1990TwodimensionalDavenportschinzelProblem}.
An example of a minimization diagram can be found in Figure~\ref{fig::mindiagramLowi}(a).
In multi-objective optimization, similar approaches have been developed.
There, the term \emph{weight set decomposition}~\cite{przybylski.gandibleux.ea2009RecursiveAlgorithmFinding} is used and the weight set is decomposed.
A formal description of the theoretical properties of minimization diagrams is given by so-called \emph{Davenport-Shinzel sequences}, see~\cite{schwartz.sharir1990TwodimensionalDavenportschinzelProblem,agarwal.sharir2000DavenportSchinzelSequences}.
Another rigorous analysis can be found in the multi-objective literature in~\cite{przybylski.gandibleux.ea2009RecursiveAlgorithmFinding}.

\medskip

Since nearly all available literature is concerned with the $2$-parametric case (e.g.,~\cite{Gusfield_1983_Parametric_Combinatorial,fernandez-baca1990SpacesweepAlgorithmsParametric,waterman.eggert.ea1992ParametricSequenceComparisons,zimmer.lengauer1997FastNumericallyStable,przybylski.gandibleux.ea2009RecursiveAlgorithmFinding,halffmann.dietz.ea2020InnerApproximationMethod}), we assume that $p=2$ for the rest of this section, unless stated otherwise.
If a problem uses more than two parameters, strategies like the ones outlined in \Cref{subsec::bensonlike} are often applied.

\medskip

Deciding if two critical regions share an edge and, if so, finding the vertices of this edge can be done by solving a single-parametric problem with any algorithm (e.g., the methods from \Cref{subsec::ESmethod,subsec::megiddogusfield}).
This problem can be seen as intersecting the parameter set $\Lambda$ with a straight line.
Hence, procedures that investigate potential edges of critical regions are called \emph{line searches}.
As notational shorthand throughout this section, we denote the number of vertices in a minimization diagram by~$V$ and let~$E$ denote the number of edges in it.

\medskip

Nearly all algorithms that compute minimization diagrams follow a similar pattern:
An initial solution (or a set of initial solutions) is computed, and its critical region is constructed by a sequence of line searches.
During these line searches, new solutions are discovered, for which the process is repeated until all critical regions have been found.

The first such algorithm in parametric literature is presented by \citeauthor{Gusfield_1983_Parametric_Combinatorial}~\cite{Gusfield_1983_Parametric_Combinatorial}.
Here, critical regions are computed by simulating ray shooting with line searches.
Rays are shot from the interior of a critical region to determine its edges.
The running time bound of the procedure is stated as $\mathcal{O}(E\cdot {T_\Pi^2})$\footnote{\citeauthor{fernandez-baca1990SpacesweepAlgorithmsParametric}~\cite{fernandez-baca1990SpacesweepAlgorithmsParametric} claims that the running time analysis by \citeauthor{Gusfield_1983_Parametric_Combinatorial}~\cite{Gusfield_1983_Parametric_Combinatorial} contains an error, but no details are given on the specific nature of this error.}.

\medskip

\citeauthor{fernandez-baca1990SpacesweepAlgorithmsParametric}~\cite{fernandez-baca1990SpacesweepAlgorithmsParametric} develops a sweep-line approach:
Edges of the minimization diagram are traced via line searches, and its vertices are enumerated in ascending order of their~$\lambda_1$ value.
The running time is in $\mathcal{O}(E\cdot {T_\Pi}^2+ E\log E)$.
A generalization for the $3$-parametric case is outlined shortly.

\medskip

\citeauthor{waterman.eggert.ea1992ParametricSequenceComparisons}~\cite{waterman.eggert.ea1992ParametricSequenceComparisons} and \citeauthor{zimmer.lengauer1997FastNumericallyStable}~\cite{zimmer.lengauer1997FastNumericallyStable} develop algorithms for genetic alignment problems, which can also be used as general strategies for finding minimization diagrams.
The first, given by \citeauthor{waterman.eggert.ea1992ParametricSequenceComparisons}~\cite{waterman.eggert.ea1992ParametricSequenceComparisons}, operates by tracing the boundary of each critical region separately.
No running time analysis is given for this algorithm.
In the algorithm from \citeauthor{zimmer.lengauer1997FastNumericallyStable}~\cite{zimmer.lengauer1997FastNumericallyStable}, a graph is used to encode adjacency relations of critical regions in the minimization diagram.
This is used to systematically check and identify potential edges.
The running time is in $\mathcal{O}(B^4+(B+V)T_\Pi)$, and a computational study in~\cite{zimmer.lengauer1997FastNumericallyStable} shows that this approach outperforms the algorithm from \citeauthor{waterman.eggert.ea1992ParametricSequenceComparisons}~\cite{waterman.eggert.ea1992ParametricSequenceComparisons} on genetic alignment problems.

\medskip

In the field of multi-objective optimization, the first algorithm that computes a minimization diagram was developed by \citeauthor{benson.sun2002WeightSetDecompositiona}~\cite{benson.sun2002WeightSetDecompositiona}, but can only be used for multi-objective linear programming problems.
The first two general algorithms were published by \citeauthor{przybylski.gandibleux.ea2009RecursiveAlgorithmFinding}~\cite{przybylski.gandibleux.ea2009RecursiveAlgorithmFinding}.
One of them is called the \emph{brute force algorithm} and works almost identically to the algorithm of \citeauthor{zimmer.lengauer1997FastNumericallyStable}~\cite{zimmer.lengauer1997FastNumericallyStable}. The other algorithm tries to improve the practical running time of the brute force algorithm by eliminating redundant line searches.
A computational study shows no significant advantage for either of the two variants.
Another algorithm is presented by \citeauthor{halffmann.dietz.ea2020InnerApproximationMethod}~\cite{halffmann.dietz.ea2020InnerApproximationMethod}.
It operates by constructing inner approximations for the critical regions and improving them iteratively through line searches.
The running time of this algorithm is in $\mathcal{O}(B^2 T_\Pi)$, and a computational study shows a similar performance to the algorithm of \citeauthor{przybylski.gandibleux.ea2009RecursiveAlgorithmFinding}~\cite{przybylski.gandibleux.ea2009RecursiveAlgorithmFinding} in a computational study.
Some further improvements are developed in~\cite{halffmann2021AdvancesMultiobjectiveOptimisation}.

\medskip

\citeauthor{alves.costa2016GraphicalExplorationWeight}~\cite{alves.costa2016GraphicalExplorationWeight} present an algorithm for computing a subset of a critical region.
The approach can be used to compute the entire minimization diagram, however, no running time bound is stated in~\cite{alves.costa2016GraphicalExplorationWeight}.

\subsubsection{Lower Image Computation}\label{subsec::bensonlike}

An approach that is closely related to parameter set decomposition is the construction of the so-called \emph{lower image}~$\lowi$.
While the concept appears earlier in the parametric literature\footnote{To the best of our knowledge, the first mention is by \citeauthor{fernandez-baca.srinivasan1991ConstructingMinimizationDiagrama}~\cite{fernandez-baca.srinivasan1991ConstructingMinimizationDiagrama}.}, the term ``lower image'' stems from multi-objective optimization~\cite{heyde.lohne2008GeometricDualityMultiple} and is the first name explicitly given to the polyhedron~$\lowi$.
Again, we directly take the parametric perspective when describing results from multi-objective optimization.

\begin{definition}[\citeauthor{heyde.lohne2008GeometricDualityMultiple}~\cite{heyde.lohne2008GeometricDualityMultiple}]
  The \emph{lower image}~$\lowi$ is the set
  \[
    \lowi\coloneqq\left\{(\lambda, z): \lambda\in\Lambda,z\leq F^*(\lambda)\right\}.
  \]
\end{definition}

The lower image can equivalently be described by the solution set:
For every solution $x\in X$, the set $H(x)=\{(\lambda,z)\in\mathbb{R}^p\times\mathbb{R}:z\leq F_\lambda(x)\}$ describes a halfspace.
The intersection of $\Lambda\times\mathbb{R}$ and all halfspaces~$H(x)$ for $x \in X$ is exactly the lower image~\cite{heyde.lohne2008GeometricDualityMultiple}.
However, the halfspaces from any optimal solution set are sufficient to obtain a minimal description of~$\lowi$.
This fact is used by several algorithmic approaches that construct~$\lowi$ using techniques from computational geometry to obtain an optimal solution set.

\medskip

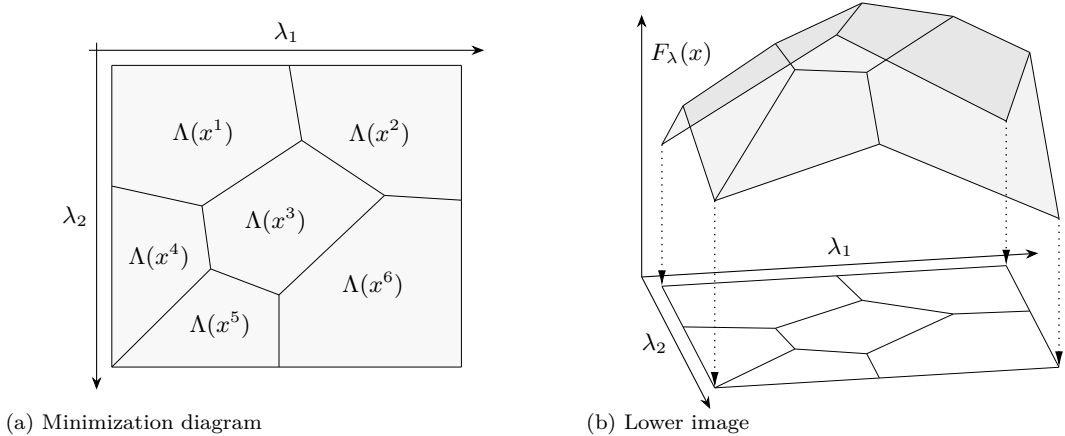
\begin{figure}
  \centering
  \begin{minipage}{0.48\textwidth}
    \centering
    \vspace{0pt}
\tdplotsetmaincoords{0}{180}
\begin{tikzpicture}[
    tdplot_main_coords,
    scale=1.000000,
    facet/.style={fill=black!100, fill opacity=0.03},
    edge/.style={color=black, line width=0.1},
    vertex/.style={},
  ]

\fill[facet] (4.000000, 0.000000) -- (4.000000, 1.600000) -- (2.813187, 1.863736) -- (1.500938, 0.994371) -- (1.666667, 0.000000) -- cycle {};
\fill[facet] (-0.600000, 0.000000) -- (1.666667, 0.000000) -- (1.500938, 0.994371) -- (0.411523, 1.724280) -- (-0.600000, 1.787500) -- cycle {};
\fill[facet] (-0.600000, 4.000000) -- (-0.600000, 1.787500) -- (0.411523, 1.724280) -- (1.800000, 3.043333) -- (1.800000, 4.000000) -- cycle {};
\fill[facet] (4.000000, 4.000000) -- (1.800000, 4.000000) -- (1.800000, 3.043333) -- (2.698795, 2.698795) -- cycle {};
\fill[facet] (4.000000, 4.000000) -- (2.698795, 2.698795) -- (2.813187, 1.863736) -- (4.000000, 1.600000) -- cycle {};
\fill[facet] (1.800000, 3.043333) -- (0.411523, 1.724280) -- (1.500938, 0.994371) -- (2.813187, 1.863736) -- (2.698795, 2.698795) -- cycle {};

\draw[edge] (1.500938, 0.994371) -- (1.666667, 0.000000);
\draw[edge] (2.813187, 1.863736) -- (4.000000, 1.600000);
\draw[edge] (2.813187, 1.863736) -- (1.500938, 0.994371);
\draw[edge] (4.000000, 1.600000) -- (4.000000, 0.000000);
\draw[edge] (4.000000, 0.000000) -- (1.666667, 0.000000);
\draw[edge] (-0.600000, 1.787500) -- (0.411523, 1.724280);
\draw[edge] (1.500938, 0.994371) -- (0.411523, 1.724280);
\draw[edge] (-0.600000, 0.000000) -- (1.666667, 0.000000);
\draw[edge] (-0.600000, 0.000000) -- (-0.600000, 1.787500);
\draw[edge] (1.800000, 4.000000) -- (1.800000, 3.043333);
\draw[edge] (0.411523, 1.724280) -- (1.800000, 3.043333);
\draw[edge] (-0.600000, 4.000000) -- (1.800000, 4.000000);
\draw[edge] (-0.600000, 1.787500) -- (-0.600000, 4.000000);
\draw[edge] (2.698795, 2.698795) -- (4.000000, 4.000000);
\draw[edge] (1.800000, 3.043333) -- (2.698795, 2.698795);
\draw[edge] (1.800000, 4.000000) -- (4.000000, 4.000000);
\draw[edge] (2.813187, 1.863736) -- (2.698795, 2.698795);
\draw[edge] (4.000000, 1.600000) -- (4.000000, 4.000000);

\node[vertex] at (2.813187, 1.863736) {};
\node[vertex] at (4.000000, 1.600000) {};
\node[vertex] at (1.500938, 0.994371) {};
\node[vertex] at (-0.600000, 0.000000) {};
\node[vertex] at (4.000000, 0.000000) {};
\node[vertex] at (1.666667, 0.000000) {};
\node[vertex] at (-0.600000, 1.787500) {};
\node[vertex] at (0.411523, 1.724280) {};
\node[vertex] at (-0.600000, 4.000000) {};
\node[vertex] at (1.800000, 4.000000) {};
\node[vertex] at (1.800000, 3.043333) {};
\node[vertex] at (2.698795, 2.698795) {};
\node[vertex] at (4.000000, 4.000000) {};


\node at ($0.2*(4.000000, 0.000000) + 0.2*(4.000000, 1.600000) + 0.2*(2.813187, 1.863736) + 0.2*(1.500938, 0.994371) + 0.2*(1.666667, 0.000000)$) {$\Lambda(x^1)$};
\node at ($0.2*(-0.600000, 0.000000) + 0.2*(1.666667, 0.000000) + 0.2*(1.500938, 0.994371) + 0.2*(0.411523, 1.724280) + 0.2*(-0.600000, 1.787500)$) {$\Lambda(x^2)$};
\node at ($0.2*(1.800000, 3.043333) + 0.2*(0.411523, 1.724280) + 0.2*(1.500938, 0.994371) + 0.2*(2.813187, 1.863736) + 0.2*(2.698795, 2.698795)$) {$\Lambda(x^3)$};
\node at ($0.25*(4.000000, 4.000000) + 0.25*(2.698795, 2.698795) + 0.25*(2.813187, 1.863736) + 0.25*(4.000000, 1.600000)$) {$\Lambda(x^4)$};
\node at ($0.25*(4.000000, 4.000000) + 0.25*(1.800000, 4.000000) + 0.25*(1.800000, 3.043333) + 0.25*(2.698795, 2.698795)$) {$\Lambda(x^5)$};
\node at ($0.2*(-0.600000, 4.000000) + 0.2*(-0.600000, 1.787500) + 0.2*(0.411523, 1.724280) + 0.2*(1.800000, 3.043333) + 0.2*(1.800000, 4.000000)$) {$\Lambda(x^6)$};

\coordinate (pseudoOrigin) at (4.2,-0.2);
\draw[-{Stealth}] ($(pseudoOrigin)+(0.1,0)$) -- ($(pseudoOrigin)+(-5.1,0)$) node[pos=0.5, above] {$\lambda_1$};
\draw[-{Stealth}] ($(pseudoOrigin)+(0,-0.1)$) -- ($(pseudoOrigin)+(0,4.5)$) node[pos=0.5,left] {$\lambda_2$};

\end{tikzpicture}
  \end{minipage}
  \hfill
  \begin{minipage}{0.48\textwidth}
    \centering
    \vspace{0pt}
  \tdplotsetmaincoords{70}{170}  \begin{tikzpicture}[
    tdplot_main_coords,    scale=1.00000,
    topfacet/.style={fill=black!100, fill opacity=0.05},
    topfocus/.style={fill=green!100, fill opacity=0.05},
    botfocus/.style={fill=green!100, fill opacity=0.05},
    botfacet/.style={fill=black!100, fill opacity=0.0},
    topedge/.style={color=black, line width=0.1},
    botedge/.style={color=black, line width=0.1},
    topvertex/.style={},
    botvertex/.style={},
    vertexarrow/.style={color=black, line width=0.5, dotted, -{Latex[width=3,length=6]}},
  ]

\fill[topfacet] (4.000000, 0.000000, 1.600000) -- (4.000000, 1.600000, 2.720000) -- (2.813187, 1.863736, 3.616703) -- (1.500938, 0.994371, 3.795497) -- (1.666667, 0.000000, 3.000000) -- cycle {};
\fill[topfacet] (-0.600000, 0.000000, 1.640000) -- (1.666667, 0.000000, 3.000000) -- (1.500938, 0.994371, 3.795497) -- (0.411523, 1.724280, 3.798765) -- (-0.600000, 1.787500, 3.248750) -- cycle {};
\fill[topfacet] (1.800000, 3.043333, 3.569667) -- (0.411523, 1.724280, 3.798765) -- (1.500938, 0.994371, 3.795497) -- (2.813187, 1.863736, 3.616703) -- (2.698795, 2.698795, 3.541205) -- cycle {};
\fill[topfacet] (4.000000, 4.000000, 2.240000) -- (2.698795, 2.698795, 3.541205) -- (2.813187, 1.863736, 3.616703) -- (4.000000, 1.600000, 2.720000) -- cycle {};
\fill[topfacet] (4.000000, 4.000000, 2.240000) -- (1.800000, 4.000000, 2.900000) -- (1.800000, 3.043333, 3.569667) -- (2.698795, 2.698795, 3.541205) -- cycle {};
\fill[topfacet] (-0.600000, 4.000000, 1.700000) -- (-0.600000, 1.787500, 3.248750) -- (0.411523, 1.724280, 3.798765) -- (1.800000, 3.043333, 3.569667) -- (1.800000, 4.000000, 2.900000) -- cycle {};


\fill[botfacet] (4.000000, 0.000000, -0.400000) -- (4.000000, 1.600000, -0.400000) -- (2.813187, 1.863736, -0.400000) -- (1.500938, 0.994371, -0.400000) -- (1.666667, 0.000000, -0.400000) -- cycle {};
\fill[botfacet] (-0.600000, 0.000000, -0.400000) -- (1.666667, 0.000000, -0.400000) -- (1.500938, 0.994371, -0.400000) -- (0.411523, 1.724280, -0.400000) -- (-0.600000, 1.787500, -0.400000) -- cycle {};
\fill[botfacet] (1.800000, 3.043333, -0.400000) -- (0.411523, 1.724280, -0.400000) -- (1.500938, 0.994371, -0.400000) -- (2.813187, 1.863736, -0.400000) -- (2.698795, 2.698795, -0.400000) -- cycle {};
\fill[botfacet] (4.000000, 4.000000, -0.400000) -- (2.698795, 2.698795, -0.400000) -- (2.813187, 1.863736, -0.400000) -- (4.000000, 1.600000, -0.400000) -- cycle {};
\fill[botfacet] (4.000000, 4.000000, -0.400000) -- (1.800000, 4.000000, -0.400000) -- (1.800000, 3.043333, -0.400000) -- (2.698795, 2.698795, -0.400000) -- cycle {};
\fill[botfacet] (-0.600000, 4.000000, -0.400000) -- (-0.600000, 1.787500, -0.400000) -- (0.411523, 1.724280, -0.400000) -- (1.800000, 3.043333, -0.400000) -- (1.800000, 4.000000, -0.400000) -- cycle {};


\draw[topedge] (1.500938, 0.994371, 3.795497) -- (1.666667, 0.000000, 3.000000);
\draw[topedge] (2.813187, 1.863736, 3.616703) -- (4.000000, 1.600000, 2.720000);
\draw[topedge] (2.813187, 1.863736, 3.616703) -- (1.500938, 0.994371, 3.795497);
\draw[topedge] (4.000000, 1.600000, 2.720000) -- (4.000000, 0.000000, 1.600000);
\draw[topedge] (4.000000, 0.000000, 1.600000) -- (1.666667, 0.000000, 3.000000);
\draw[topedge] (-0.600000, 1.787500, 3.248750) -- (0.411523, 1.724280, 3.798765);
\draw[topedge] (1.500938, 0.994371, 3.795497) -- (0.411523, 1.724280, 3.798765);
\draw[topedge] (-0.600000, 0.000000, 1.640000) -- (1.666667, 0.000000, 3.000000);
\draw[topedge] (-0.600000, 0.000000, 1.640000) -- (-0.600000, 1.787500, 3.248750);
\draw[topedge] (1.800000, 4.000000, 2.900000) -- (1.800000, 3.043333, 3.569667);
\draw[topedge] (0.411523, 1.724280, 3.798765) -- (1.800000, 3.043333, 3.569667);
\draw[topedge] (-0.600000, 4.000000, 1.700000) -- (1.800000, 4.000000, 2.900000);
\draw[topedge] (-0.600000, 1.787500, 3.248750) -- (-0.600000, 4.000000, 1.700000);
\draw[topedge] (2.698795, 2.698795, 3.541205) -- (4.000000, 4.000000, 2.240000);
\draw[topedge] (1.800000, 3.043333, 3.569667) -- (2.698795, 2.698795, 3.541205);
\draw[topedge] (1.800000, 4.000000, 2.900000) -- (4.000000, 4.000000, 2.240000);
\draw[topedge] (2.813187, 1.863736, 3.616703) -- (2.698795, 2.698795, 3.541205);
\draw[topedge] (4.000000, 1.600000, 2.720000) -- (4.000000, 4.000000, 2.240000);

\draw[botedge] (1.500938, 0.994371, -0.400000) -- (1.666667, 0.000000, -0.400000);
\draw[botedge] (2.813187, 1.863736, -0.400000) -- (4.000000, 1.600000, -0.400000);
\draw[botedge] (2.813187, 1.863736, -0.400000) -- (1.500938, 0.994371, -0.400000);
\draw[botedge] (4.000000, 1.600000, -0.400000) -- (4.000000, 0.000000, -0.400000);
\draw[botedge] (4.000000, 0.000000, -0.400000) -- (1.666667, 0.000000, -0.400000);
\draw[botedge] (-0.600000, 1.787500, -0.400000) -- (0.411523, 1.724280, -0.400000);
\draw[botedge] (1.500938, 0.994371, -0.400000) -- (0.411523, 1.724280, -0.400000);
\draw[botedge] (-0.600000, 0.000000, -0.400000) -- (1.666667, 0.000000, -0.400000);
\draw[botedge] (-0.600000, 0.000000, -0.400000) -- (-0.600000, 1.787500, -0.400000);
\draw[botedge] (1.800000, 4.000000, -0.400000) -- (1.800000, 3.043333, -0.400000);
\draw[botedge] (0.411523, 1.724280, -0.400000) -- (1.800000, 3.043333, -0.400000);
\draw[botedge] (-0.600000, 4.000000, -0.400000) -- (1.800000, 4.000000, -0.400000);
\draw[botedge] (-0.600000, 1.787500, -0.400000) -- (-0.600000, 4.000000, -0.400000);
\draw[botedge] (2.698795, 2.698795, -0.400000) -- (4.000000, 4.000000, -0.400000);
\draw[botedge] (1.800000, 3.043333, -0.400000) -- (2.698795, 2.698795, -0.400000);
\draw[botedge] (1.800000, 4.000000, -0.400000) -- (4.000000, 4.000000, -0.400000);
\draw[botedge] (2.813187, 1.863736, -0.400000) -- (2.698795, 2.698795, -0.400000);
\draw[botedge] (4.000000, 1.600000, -0.400000) -- (4.000000, 4.000000, -0.400000);

\node[topvertex] at (2.813187, 1.863736, 3.616703) {};
\node[topvertex] at (4.000000, 1.600000, 2.720000) {};
\node[topvertex] at (1.500938, 0.994371, 3.795497) {};
\node[topvertex] at (-0.600000, 0.000000, 1.640000) {};
\node[topvertex] at (4.000000, 0.000000, 1.600000) {};
\node[topvertex] at (1.666667, 0.000000, 3.000000) {};
\node[topvertex] at (-0.600000, 1.787500, 3.248750) {};
\node[topvertex] at (0.411523, 1.724280, 3.798765) {};
\node[topvertex] at (-0.600000, 4.000000, 1.700000) {};
\node[topvertex] at (1.800000, 4.000000, 2.900000) {};
\node[topvertex] at (1.800000, 3.043333, 3.569667) {};
\node[topvertex] at (2.698795, 2.698795, 3.541205) {};
\node[topvertex] at (4.000000, 4.000000, 2.240000) {};

\node[botvertex] at (2.813187, 1.863736, -0.400000) {};
\node[botvertex] at (4.000000, 1.600000, -0.400000) {};
\node[botvertex] at (1.500938, 0.994371, -0.400000) {};
\node[botvertex] at (-0.600000, 0.000000, -0.400000) {};
\node[botvertex] at (4.000000, 0.000000, -0.400000) {};
\node[botvertex] at (1.666667, 0.000000, -0.400000) {};
\node[botvertex] at (-0.600000, 1.787500, -0.400000) {};
\node[botvertex] at (0.411523, 1.724280, -0.400000) {};
\node[botvertex] at (-0.600000, 4.000000, -0.400000) {};
\node[botvertex] at (1.800000, 4.000000, -0.400000) {};
\node[botvertex] at (1.800000, 3.043333, -0.400000) {};
\node[botvertex] at (2.698795, 2.698795, -0.400000) {};
\node[botvertex] at (4.000000, 4.000000, -0.400000) {};

\draw[vertexarrow] (4.000000, 0.000000, 1.600000) -- (4.000000, 0.000000, -0.400000);
\draw[vertexarrow] (4.000000, 4.000000, 2.240000) -- (4.000000, 4.000000, -0.400000);
\draw[vertexarrow] (-0.600000, 4.000000, 1.700000) -- (-0.600000, 4.000000, -0.400000);
\draw[vertexarrow] (-0.600000, 0.000000, 1.640000) -- (-0.600000, 0.000000, -0.400000);


\coordinate (pseudoOrigin) at (4.2,-0.4,-0.4);
\draw[-{Stealth}] (pseudoOrigin) -- ($(pseudoOrigin)+(-5.3,0,0)$) node[pos=0.5, above] {$\lambda_1$};
\draw[-{Stealth}] (pseudoOrigin) -- ($(pseudoOrigin)+(0,5.2,0)$) node[pos=0.5,left] {$\lambda_2$};
\draw[-{Stealth}] (pseudoOrigin) -- ($(pseudoOrigin)+(0,0,3.7)$) node[pos=0.85,right] {$F_\lambda(x)$};


\end{tikzpicture}
  \end{minipage}
  \begin{minipage}[t]{0.48\textwidth}
    \subcaption{Minimization diagram}\label{fig::mindiagram}
  \end{minipage}
  \hfill
  \begin{minipage}[t]{0.48\textwidth}
    \subcaption{Lower image}\label{fig::lowi}
  \end{minipage}
  \caption{The minimization diagram and the lower image for an instance with $B=6$.
  The shadow of the lower image is exactly the minimization diagram.}\label{fig::mindiagramLowi}
\end{figure}

The concept of the minimization diagrams and the lower image are closely related.
A minimization diagram can be obtained by projecting~$\lowi$ onto the $\lambda$ plane, and~$\lowi$ can be constructed by ``lifting'' each critical region in the diagram with its objective value.
This is illustrated in \Cref{fig::mindiagramLowi}.

\medskip

Approaches that construct~$\lowi$ all follow a similar pattern:
Solutions are found iteratively, and, in each iteration, the currently known solutions are used to construct a polyhedral approximation of~$\lowi$.
From this approximation, parameter values are obtained, and the corresponding non-parametric problem is solved for these values to either find new solutions or determine that an optimal solution set has been found.
The first such algorithm is limited to the $2$-parametric case and is described by \citeauthor{fernandez-baca.srinivasan1991ConstructingMinimizationDiagrama}~\cite{fernandez-baca.srinivasan1991ConstructingMinimizationDiagrama}.
Here, the polyhedral approximation is maintained by a graph structure.
The algorithm is generalized by \citeauthor{fernandez-baca.seppalainen.ea2004ParametricMultipleSequence}~\cite{fernandez-baca.seppalainen.ea2004ParametricMultipleSequence} to the general $p$-parametric case\footnote{
At this point in time, the algorithm \citeauthor{fernandez-baca.seppalainen.ea2004ParametricMultipleSequence}~\cite{fernandez-baca.seppalainen.ea2004ParametricMultipleSequence} describe is already considered as ``part of the folklore''.
But to the best of our knowledge, no formal description of the algorithm appears in literature before~\cite{fernandez-baca.seppalainen.ea2004ParametricMultipleSequence}
} by using a procedure described by \citeauthor{dobkin.edelsbrunner.ea1986ProbingConvexPolytopes}~\cite{dobkin.edelsbrunner.ea1986ProbingConvexPolytopes} for the polyhedral approximation.

\medskip

In the multi-objective optimization literature, the so-called \emph{dual-Benson algorithm} is used to construct~$\lowi$.
First presented by \citeauthor{ehrgott.lohne.ea2012DualVariantBenson}~\cite{ehrgott.lohne.ea2012DualVariantBenson}, it was further developed in~\cite{hamel.lohne.ea2014BensonTypeAlgorithms,lohne.weissing2017VectorLinearProgram}.
Originally designed for multi-objective linear programs, it was subsequently generalized to apply to combinatorial problems~\cite{bokler.mutzel2015OutputSensitiveAlgorithmsEnumerating}, integer linear programs~\cite{borndorfer.schenker.ea2016PolySCIP}, and mixed integer linear programs~\cite{bokler.nemesch.ea2023PaMILOSolverMultiobjective}.
In the dual-Benson algorithm, the well-known \emph{vertex enumeration problem}~\cite{chazelle1993OptimalConvexHull} is used for the polyhedral approximation.
\citeauthor{bokler.mutzel2015OutputSensitiveAlgorithmsEnumerating}~\cite{bokler.mutzel2015OutputSensitiveAlgorithmsEnumerating} show that the dual-Benson algorithm runs in output-polynomial time as long as $\Pi(\lambda)$ can be solved in polynomial time for every~$\lambda\in\Lambda$.

\medskip

In multi-objective optimization, there are also a number of algorithms that compute the so-called \emph{Edgeworth-Pareto} hull (often also referred to as the \emph{upper image}).
The Edgeworth-Pareto hull is the polyhedron $\operatorname{conv}(Y)+\mathbb{R}_{\geq0}^d$, which is dual to~$\lowi$.
This duality relation is described in detail by \citeauthor{heyde.lohne2008GeometricDualityMultiple}~\cite{heyde.lohne2008GeometricDualityMultiple} for linear programming, but also applies to general mixed integer linear programming~\cite{bokler.parragh.ea2024OuterApproximationAlgorithm}.
Because of this duality relation, every algorithm that computes the Edgeworth-Pareto hull can also be used to compute~$\lowi$.
Such algorithms can be found in~\cite{benson1998OuterApproximationAlgorithm,ozpeynirci.koksalan2010ExactAlgorithmFinding,przybylski.klamroth.ea2019SimpleEfficientDichotomic,csirmaz2021InnerApproximationAlgorithm,bokler.parragh.ea2024OuterApproximationAlgorithm}.

\subsection{General Approximation Algorithms for Parametric Optimization Problems}\label{subsec::generalAPX}

General approximation algorithms for parametric problems usually work under rather mild assumptions (see, for example,~\cite{Bazgan_2022_An_approximation,Helfrich_2022_An_approximation}):

\medskip

\begin{enumerate}
  \item For every $i\in[p]$, the function~$f_i$ is non-negative and $f$ is computable in polynomial time.
  Furthermore, $F^*(\lambda)\geq0\, \forall \lambda\in\Lambda$. 
  \item An approximation oracle~$\oracle_{\beta}$ with polynomial running time is available.
  \item Often, the structure of~$\Lambda$ is further constrained, e.g., by setting $\Lambda=\mathbb{R}^p_{\geq0}$~\cite{helfrich.ruzika.ea2024EfficientlyConstructing} or by bounding each parameter value from below~\cite{Bazgan_2022_An_approximation,Helfrich_2022_An_approximation}.
\end{enumerate}

\medskip

There exist only a few general approximation algorithms for parametric optimization problems.
For the single-parametric case, the adaption of the Eisner-Severance method by \citeauthor{Bazgan_2022_An_approximation}~\cite{Bazgan_2022_An_approximation} is available (see \Cref{subsec::ESmethod}).

\medskip

For multi-parametric problems, an early approach is described by \citeauthor{katoh.ibaraki1987ParametricCharacterizationEapproximation}~\cite{katoh.ibaraki1987ParametricCharacterizationEapproximation}, but it is only used as part of a larger algorithm there.
A dedicated approximation algorithm for multi-parametric problems is developed by \citeauthor{Helfrich_2022_An_approximation}~\cite{Helfrich_2022_An_approximation}.
The approaches of both \citeauthor{katoh.ibaraki1987ParametricCharacterizationEapproximation}~\cite{katoh.ibaraki1987ParametricCharacterizationEapproximation} and \citeauthor{Helfrich_2022_An_approximation}~\cite{Helfrich_2022_An_approximation} work by constructing a grid over the parameter set.
For each vertex of this grid, the non-parametric problem is solved using~$\oracle_\beta$.
The structure of the grid then ensures that the returned solution is $(\alpha\cdot\beta)$-approximate for some $\alpha>1$ for parameter values in a small neighborhood of this vertex.
Furthermore, the union of these neighborhoods covers the entire parameter set, so the set of all solutions returned by the oracle calls at the vertices is an $(\alpha\cdot\beta)$-approximation set.
If the oracle is exact or an (F)PTAS, the parametric approximation algorithm itself is an (F)PTAS.
When comparing the algorithm of \citeauthor{katoh.ibaraki1987ParametricCharacterizationEapproximation}~\cite{katoh.ibaraki1987ParametricCharacterizationEapproximation} to that of \citeauthor{helfrich.ruzika.ea2024EfficientlyConstructing}~\cite{helfrich.ruzika.ea2024EfficientlyConstructing}, the latter needs fewer oracle calls and offers more flexibility regarding the choice of the parameter set.

\medskip

Some results on the cardinality of $\alpha$-approximation sets are presented by \citeauthor{Helfrich_2022_An_approximation}~\cite{Helfrich_2022_An_approximation}.
First, even when given access to an exact oracle $\oracle_{1}$, no algorithm that has problem access only through this oracle can approximate the optimal cardinality~$B$ better than by a factor of~$p+1$.
For $p=1$, \citeauthor{Bazgan_2022_An_approximation}~\cite{Bazgan_2022_An_approximation} describe an algorithm that can achieve the best-possible factor of~$2$.
Second, for $p=3$ and a routine $\oracle_{\beta}$ with $\beta>1$, the cardinality $B$ cannot be approximated within any constant factor~$L\in\mathbb{N}$ by an algorithm that only has access to this routine.

\medskip

As mentioned in \Cref{subsec::moopt}, a survey of approximation algorithms for multi-objective optimization problems is given by \citeauthor{herzel.ruzika.ea2021ApproximationMethodsMultiobjective}~\cite{herzel.ruzika.ea2021ApproximationMethodsMultiobjective}.
Every $\alpha$-convex approximation algorithm mentioned therein can also be used for parametric optimization.
A recent algorithm for $\alpha$-convex approximation is developed by \citeauthor{helfrich.ruzika.ea2024EfficientlyConstructing}~\cite{helfrich.ruzika.ea2024EfficientlyConstructing}.
It builds upon the dual-Benson algorithm from multi-objective optimization.
In contrast to the grid-based approaches, it chooses the parameter values where $\oracle_{\beta}$ is called in an adaptive fashion.
A computational study in~\cite{helfrich.ruzika.ea2024EfficientlyConstructing} shows that this adaptive strategy outperforms other general approximation algorithms on a set of generated instances for parametric knapsack and metric traveling salesman problems.

\subsection{Parametric Mixed Integer Linear Programming}\label{subsec::milp}

A \emph{mixed integer linear program} (MILP) is the problem
\begin{equation}\label{milp::formulation}
\begin{array}{l r c r c l}
    \mathrm{min} & c^\T x&  + & d^\T y \\
    \mathrm{s.t.} & A x & + & B y & \geq & b, \\
    & \multicolumn{5}{l}{x\in\mathbb{R}^{n_{\mathbb{R}}},y\in\mathbb{Z}^{n_{\mathbb{Z}}},}
\end{array}
\end{equation}

where $c\in\mathbb{R}^{n_{\mathbb{R}}}$, $d\in\mathbb{R}^{n_{\mathbb{Z}}}$, $A\in\mathbb{R}^{m_{\mathbb{R}}\times n_{\mathbb{R}}}$, $B\in\mathbb{R}^{m_{\mathbb{Z}}\times n_{\mathbb{Z}}}$, and $b\in\mathbb{R}^{m_{\mathbb{R}}+m_{\mathbb{Z}}}$ for some $m_{\mathbb{R}},m_{\mathbb{Z}},n_{\mathbb{R}},n_{\mathbb{Z}}\geq0$. If $n_{\mathbb{Z}}=0$, the problem is a \emph{linear program} (LP), and if $n_{\mathbb{R}}=0$, it is an \emph{integer linear program} (ILP).
Parameter dependencies may appear in the objective function, the constraint matrix, or the right-hand side vector.

\medskip

Note that all of the specific problems discussed in \Cref{subsec::specificProblems} can also be formulated as a parametric MILP.
In this section, we take a short look at the case where a problem is explicitly formulated as in~\eqref{milp::formulation}.
As this is an extensively researched area of parametric optimization, a complete coverage is out of scope for this survey.
Instead, we point to some selected surveys, which together provide a comprehensive overview.

For parametric LPs, an introduction to several parametric simplex methods is provided by \citeauthor{gal1995PostoptimalAnalysesParametric}~\cite{gal1995PostoptimalAnalysesParametric}.
For MILPs, two surveys covering early developments are provided by \citeauthor{Geoffrion_1977_Exceptional_Paper}~\cite{Geoffrion_1977_Exceptional_Paper} and \citeauthor{greenberg1998AnnotatedBibliographyPostSolution}~\cite{greenberg1998AnnotatedBibliographyPostSolution}.
Two surveys on recent developments are presented by \citeauthor{Pistikopoulos_2012_Theoretical_and}~\cite{Pistikopoulos_2012_Theoretical_and} and \citeauthor{oberdieck.diangelakis.ea2016MultiparametricProgrammingIts}~\cite{oberdieck.diangelakis.ea2016MultiparametricProgrammingIts}.
They summarize the relevant literature and provide some further introduction into key techniques for single- and multi-parametric LPs, ILP, and MILP (and some more specialized problem classes).

\medskip

An important result on the cardinality of optimal solution sets for parametric LPs is given by \citeauthor{Murty_1980_Computational_Complexity}~\cite{Murty_1980_Computational_Complexity}.
A single-parametric LP with a parameter-dependent right-hand side vector is constructed, where~$n_{\mathbb{R}}$ is even, $m_{\mathbb{R}}=\nicefrac{n_{\mathbb{R}}}{2}$, and the cardinality of an optimal solution set is $2^{\nicefrac{n_{\mathbb{R}}}{2}}$.
Due to LP duality, this result also holds for parametric LPs where only the objective function is parameter-dependent.
Furthermore, since every parametric LP is also a parametric MILP, the same lower bound applies to parametric MILPs in general.


\section{Specific Problems}\label{subsec::specificProblems}

In the following sections, we take a look at several parametric combinatorial optimization problems that are well-researched in the literature.
We discuss the shortest path problem in Section~\ref{subsec::ShortestPath}, the minimum cut (or maximum flow) problem in Section~\ref{subsec::MinCut}, the minimum cost flow problem in Section~\ref{subsec::MinCostFlow}, the matroid problem and minimum spanning tree problem in Section~\ref{subsec::MinSpannTree}, the knapsack problem in Section~\ref{subsec::Knapsack}, and the matching and assignment problems in Section~\ref{subsec::Matching}.
We list bounds on the cardinality of~$B$ and provide an overview of important algorithms for each problem.

\subsection{Parametric Shortest Path Problems}
\label{subsec::ShortestPath}

A parametric shortest path problem consists of a directed connected graph~$G=(V,A)$ with $|V|=n$ vertices, $|A|=m$ arcs, a parameter set $\Lambda \subseteq \mathbb{R}^p$, and parametric arc weight functions $w_{a}(\lambda)=w_0(a) + \sum_{i=1}^{p}\lambda_i \cdot w_i(a)$ for~$\lambda\in\Lambda$ and~$a \in A$.
There are three well-known variants of parametric shortest path problems.

The first variant is the parametric \emph{single-pair shortest path problem} (SPSP), where a source vertex~$s\in V$ and a target vertex~$t\in V$ are given.
A feasible solution is a directed simple path from~$s$ to~$t$, and the goal is to compute, for every $\lambda\in\Lambda$, a feasible $s$-$t$-path~$P$ with minimum total length $F_\lambda(x)=\sum_{a \in P}w_a(\lambda)$.

The second variant is the parametric \emph{single-source shortest path problem} (SSSP), where only a source vertex~$s\in V$ is given.
A feasible solution is a set of paths from~$s$ to all other vertices in~$V$, and the goal is to find, for every $\lambda\in\Lambda$, a set of paths such that each of them is a shortest path to its target vertex.

The third variant is the parametric \emph{all-pairs shortest path problem} (APSP), where a feasible solution is a set of paths between all pairs of vertices. The goal is to find, for every $\lambda\in\Lambda$, a set of paths such that all pairs of vertices are connected by a shortest path.

\subsubsection{Intractability Results}\label{sec::shortestpath::intractability}

Since even the cardinality of optimal solution sets for the single-parametric SPSP can be exponential in the size of the graph, parametric shortest path problems are intractable in general.
This result is shown in several papers.

The first proof is given by \citeauthor{Carstensen_1983_Complexity_of}~\cite{Carstensen_1983_Complexity_of}.
She proves that, even in the single-parametric setting, there are instances for which the cardinality of an optimal solution set is in~$n^{\Omega(\log(n))}$.
\citeauthor{Carstensen_1983_The_Complexity}~\cite{Carstensen_1983_The_Complexity} also proves that this lower bound is tight.
A matching upper bound of~$n^{\log(n)+\mathcal{O}(1)}$ on the cardinality of an optimal solution set is attributed to \citeauthor{Gusfield_1980_Sensitivity_Analysis}~\cite{Gusfield_1980_Sensitivity_Analysis}.

\medskip

\citeauthor{Mulmuley_2001_A_Lower}~\cite{Mulmuley_2001_A_Lower} simplify the proof from \citeauthor{Carstensen_1983_Complexity_of}~\cite{Carstensen_1983_Complexity_of} and extend it by taking into account the bit-lengths of the arc weight functions. There is a revised variant of this proof \footnote{\citeauthor{Gajjar_2019_Parametric_Shortest}~\cite{Gajjar_2019_Parametric_Shortest} remark that the proof from \citeauthor{Mulmuley_2001_A_Lower}~\cite{Mulmuley_2001_A_Lower} ``is rather cryptic and has errors that throw the reader off''. } from \citeauthor{Gajjar_2019_Parametric_Shortest}~\cite{Gajjar_2019_Parametric_Shortest}.
Their result transfers the lower bound of $n^{\Omega(\log(n))}$ also to planar graphs.
Thereby, they refute an earlier conjecture by \citeauthor{Nikolova_2009_Strategic_Algorithms}~\cite{Nikolova_2009_Strategic_Algorithms}.
The proof idea of \citeauthor{Gajjar_2019_Parametric_Shortest}~\cite{Gajjar_2019_Parametric_Shortest} is to inductively compose dense bipartite graphs, so that the number of vertices increases by a constant factor and the number of breakpoints of the optimal cost function increases by a linear factor.
Meanwhile, these non-planar bipartite graphs can be simulated by a planar gadget, where each arc is replaced by a path consisting of up to $n^2$~arcs and the original weight is distributed among these arcs.

\subsubsection{Tractability Results for Special Cases}

While parametric shortest path problems are, in general, intractable, some special cases allow polynomial-time algorithms.

\citeauthor{Karp_1981_Parametric_shortest}~\cite{Karp_1981_Parametric_shortest} consider single-parametric SSSPs in which the parameter coefficient~$w_1(a)$ is either~$-1$ or~$0$ for each arc~$a\in A$.
The parameter set~$\Lambda$ is chosen such that there is no negative-cost cycle for any~$\lambda\in\Lambda$.
Two algorithms are provided in~\cite{Karp_1981_Parametric_shortest}, one with a running time in $\mathcal{O}\big(n^3\big)$, the other in $\mathcal{O}\big(n m \log(n)\big)$.
The latter algorithm works by maintaining a shortest path tree that changes iteratively while increasing the parameter.
\citeauthor{Young_1991_Faster_parametric}~\cite{Young_1991_Faster_parametric} improve the running time bound of this algorithm by using Fibonacci heaps~\cite{Fredman_1987_Fibonacci_heaps} to $\mathcal{O}\big(n m + n^2\cdot\log(n)\big)$.

\medskip

An algorithm for the multi-parametric generalization of the problem considered by \citeauthor{Karp_1981_Parametric_shortest}~\cite{Karp_1981_Parametric_shortest} is provided by \citeauthor{Seree_2021_An_algorithm}~\cite{Seree_2021_An_algorithm}.
The running time of this algorithm, however, is polynomial in the number of arcs, but exponential in the number of vertices.

\medskip

The parametric shortest path algorithms from \citeauthor{Karp_1981_Parametric_shortest}~\cite{Karp_1981_Parametric_shortest} and \citeauthor{Young_1991_Faster_parametric}~\cite{Young_1991_Faster_parametric} are used by \citeauthor{Erickson_2010_Maximum_Flows}~\cite{Erickson_2010_Maximum_Flows} to solve the non-parametric maximum flow problem in planar graphs.
The maximum flow problem is interpreted as a single-parametric shortest path problem in the dual graph with parameter coefficients~$w_1(a)\in\{-1,0,1\}$ for all arcs~$a\in A$.
\citeauthor{Erickson_2010_Maximum_Flows}~\cite{Erickson_2010_Maximum_Flows} provides an algorithm that solves this parametric problem in $\mathcal{O}\big(n \log(n)\big)$.
This running time bound does not hold if the algorithm is used for more general parametric shortest path problems~\cite{Erickson_2010_Maximum_Flows}.

\medskip

\citeauthor{Correa_2017_Fare_Evasion}~\cite{Correa_2017_Fare_Evasion} tackle the issue of fare evasion in transit networks by building on results for the parametric shortest path problem.
For single-parametric shortest path problems, they show that there are at most $m-1$ breakpoints on series-parallel graphs.

\subsubsection{Results for Related Shortest Path Problems}

Apart from \citeauthor{Correa_2017_Fare_Evasion}~\cite{Correa_2017_Fare_Evasion}, there are also other optimization papers that use the connection to parametric shortest path problems to achieve new results for related problems. \medskip

\citeauthor{Nikolova_2006_Stochastic_Shortest}~\cite{Nikolova_2006_Stochastic_Shortest} give an exact $n^{\Theta(\log n)}$ algorithm for the normally distributed stochastic shortest path problem, based on quasi-convex maximization. They use the connection between computational geometry and parametric optimization.
\medskip

\citeauthor{Chakraborty_2010_Two-phase_algorithms}~\cite{Chakraborty_2010_Two-phase_algorithms} construct 
two-phase algorithms for single-parametric SSSPs.
In the first phase of such an algorithm, a special data structure is created.
In the second phase, this structure is used  as an oracle that takes any $\lambda\in\Lambda$ for which there is no negative cycle as input and returns an optimal solution of $\Pi(\lambda)$.
For single-parametric SSSPs, an algorithm is presented where the first phase runs in $\Tilde{\mathcal{O}}(n^{4})$ and queries in the second phase can be answered in $\mathcal{O}(m+n\log n)$.
Another algorithm for $p$-parametric SPSPs achieves a running time in~$\mathcal{O}\left(n^{\log n}\right)$ for the first phase and in~$\mathcal{O}(\log^2 n)$ for every second phase query.
The latter algorithm can be generalized for arc weights that are polynomial functions of $\lambda$.

\medskip

\citeauthor{Foschini_2014_On_the}~\cite{Foschini_2014_On_the} analyze the complexity of time-dependent shortest path problems by utilizing connections to parametric shortest path problems.

\subsubsection{Approximation Methods}

Concerning the approximation of parametric shortest path problems, we are only aware of literature from multi-objective optimization.

For the bi-objective shortest path problem, \citeauthor{Diakonikolas_2011_Approximation_of}~\cite{Diakonikolas_2011_Approximation_of} gives a lower bound on the cardinality of approximation sets that can be computed in polynomial time.
Transferred to the parametric setting, the theorem states that there is no polynomial-time $(1+\varepsilon)$-approximation algorithm that can approximate the cardinality of the smallest $(1+\varepsilon)$-approximation set to a factor smaller than~$2$, unless $\mathsf{P}=\mathsf{NP}$.
An algorithm that achieves the best-possible factor of~$2$ can be obtained from the general approximation algorithm given by \citeauthor{Bazgan_2022_An_approximation}~\cite{Bazgan_2022_An_approximation}.

\medskip

For an overview of approximation algorithms for multi-objective shortest path problems, we refer to the recent paper by \citeauthor{Bökler_2010_Approximating_Multiobjective}~\cite{Bökler_2010_Approximating_Multiobjective} and the general survey by \citeauthor{herzel.ruzika.ea2021ApproximationMethodsMultiobjective}~\cite{herzel.ruzika.ea2021ApproximationMethodsMultiobjective}.

\subsection{Parametric Minimum Cut (Maximum Flow) Problems}\label{subsec::MinCut}

A \emph{parametric minimum cut problem} (sometimes also \emph{parametric minimum $s$-$t$-cut problem}) consists of a directed graph~$G=(V,A)$ with $|V|=n$ vertices and $|A|=m$ arcs, a source vertex~$s\in V$, a target vertex~$t\in V$ with $s \neq t$, a parameter set $\Lambda \subseteq \mathbb{R}^p$, and parametric arc capacity functions $c_{a}(\lambda)=c_0(a) + \sum_{i=1}^{p}\lambda_i \cdot c_i(a)$ for~$\lambda\in\Lambda$ and~$a \in A$.
A feasible solution~$x\subseteq \{0,1\}^m$ represents a set of arcs that partitions~$V$ into two disjointed subsets~$S$ with~$s\in S$ and~$V\setminus S$ with~$t\in V\setminus S$, a so-called \emph{$s$-$t$-cut}.
The goal is to find, for every~$\lambda \in \Lambda$, an $s$-$t$-cut~$x\subseteq{A}$ with minimum total arc capacity~$C_{\lambda}(x)=\sum_{a\in x}c_a(\lambda)$.

Due to the well-known max-flow-min-cut theorem (see~\cite{Ahuja_1993_Network_Flows}), intractability results for the parametric minimum cut problem can be applied to the parametric maximum flow problem, and vice versa.

\subsubsection{Intractability Results and Approximation Methods}

The first proof of intractability for the single-parametric minimum cut problem is given by \citeauthor{Carstensen_1983_Complexity_of}~\cite{Carstensen_1983_Complexity_of} (see \Cref{sec::shortestpath::intractability}).
However, the instance constructed in~\cite{Carstensen_1983_Complexity_of} includes negative arc capacities.
When negative arc capacities are allowed, the minimum cut problem is NP-hard.
Therefore, the proof of intractability by \citeauthor{Mulmuley_1999_Lower_bounds}~\cite{Mulmuley_1999_Lower_bounds} provides a stronger result, as it applies to parametric minimum cut problems with non-negative arc capacities in at least one parameter interval.
Any problem with a fixed parameter value in this interval can be solved in polynomial time, but there can still be superpolynomially many breakpoints within the interval.
\citeauthor{Gajjar_2019_Parametric_Shortest}~\cite{Gajjar_2019_Parametric_Shortest} tighten this result further:
They construct the planar dual of the graph used for proving intractability of the shortest path problem.
This instance can be used to verify that, for every~$k \in \mathbb{N}$, there is an instance of the single-parametric minimum cut problem consisting of non-negative integer arc capacities and a planar graph with a number of vertices polynomial in~$k$, bounded bit lengths of the arc capacity coefficients, and a linear capacity assignment to the arcs such that there are $k^{\Omega(\log(k))}$ breakpoints in the optimal cost curve.

\medskip

Using a graph from \citeauthor{Zadeh_1973_A_bad}~\cite{Zadeh_1973_A_bad}, \citeauthor{Ruhe_1988_Complexity_results}~\cite{Ruhe_1988_Complexity_results} proves intractability also for the generalized parametric maximum flow problem where the flow along each arc~$(i,j) \in A$ is scaled by a multiplier~$g_{ij} > 0$. Here, the goal is to find, for every~$\lambda \in \Lambda$, a generalized $s$-$t$-flow such that the amount of flow reaching the target vertex~$t$ is maximized.

\medskip

\citeauthor{Allman_2022_Complexity_of}~\cite{Allman_2022_Complexity_of} show that $2$-parametric minimum cut problems satisfying the so-called \emph{specialized source-sink monotonicity (sSSM)} can still have exponentially many breakpoints.
A graph with parametric arc capacities satisfies sSSM if, for each arc~$(s,j)$, the arc capacity~$c_{(s,j)}(\lambda_1,\lambda_2)$ is increasing in~$\lambda_1$ and~$\lambda_2$, and, for all other arcs, the arc capacity is independent of~$\lambda_1$ and~$\lambda_2$.

\subsubsection{Exact Solution Methods for Special Cases}

The parametric minimum cut problem is one of the most-researched problems in the parametric optimization literature. A wide range of applications is studied (e.g.,~\cite{Kolmogorov_2007_Applications_of, Asad_2013_Implementing_a, Sangeorzan_2016_A_preflow}).
Reviewing all the applications goes beyond the scope of this survey.

Instead, we focus on important exact solution methods for special cases of the problem. 
These special cases are based on the following properties:

\begin{definition}[Nesting Property]\label{Def:NestingProp}
    An instance $(c,(G,s,t),\Lambda)$ has the \emph{nesting property} if, for every sequence $\lambda_1 < \lambda_2 < \dots < \lambda_h$ of parameter values from $\Lambda$ and $h \in \mathbb{N}$, there exist corresponding (\emph{nested}) minimum cuts $(S_1,V\setminus S_1), (S_2,V\setminus S_2), \dots, (S_h,V\setminus S_h)$ such that $S_1 \subseteq S_2 \subseteq \dots \subseteq S_h$.
\end{definition}

\begin{definition}[Generalized Nesting Property]\label{Def:GenNestingProp}
    An instance $(c,(G,s,t),\Lambda)$ has the \emph{generalized nesting property} if, for every sequence $\lambda_1 < \lambda_2 < \dots < \lambda_h$  of parameter values from $\Lambda$ and $h \in \mathbb{N}$, there exist corresponding minimum cuts $(S_1,V\setminus S_1), (S_2,V\setminus S_2), \dots, (S_h,V\setminus S_h)$ such that $S_{\sigma(1)} \subseteq S_{\sigma(2)} \subseteq \dots \subseteq S_{\sigma(h)}$ for some permutation $\sigma$ on $[h]$.
\end{definition}

\begin{definition}[Source-Sink Monotonicity (SSM)]\label{Def:SSM}
    An instance of the single-parametric minimum cut problem is \emph{source-sink monotone} if, for every arc~$(i,j)$, the arc capacity~$c_{(i,j)}(\lambda)$ satisfies the condition
    \begin{equation*}
        \begin{aligned}
             c_{(i,j)}(\lambda) \text{ is } \left\{ \begin{array}{rcl}
                \text{non-decreasing in $\lambda$} &\text{ if } i=s,\,j\neq t, \\
                \text{non-increasing in $\lambda$} &\text{ if } j=t,\,i\neq s, \\
                \text{independent of $\lambda$} &\text{ otherwise.}
            \end{array}\right.
        \end{aligned}
    \end{equation*}
\end{definition}

\citeauthor{Eisner_1976_Mathematical_Techniques}~\cite{Eisner_1976_Mathematical_Techniques} and \citeauthor{Stone_1978_Critical_Load}~\cite{Stone_1978_Critical_Load} show that minimum cuts of the single-parametric minimum cut problem satisfying SSM have the nesting property.
The (generalized) nesting property implies that the optimal cost curve has at most $n-1$ breakpoints~\cite{Eisner_1976_Mathematical_Techniques}.
For instances satisfying SSM, \citeauthor{Gallo_1989_A_Fast}~\cite{Gallo_1989_A_Fast} formulate a procedure that extends the push-relabel algorithm of \citeauthor{Goldberg_1988_A_new}~\cite{Goldberg_1988_A_new} and determines these breakpoints in the same asymptotic running time as a single maximum flow computation. 

\medskip

\citeauthor{Arai_1993_Generalization_of}~\cite{Arai_1993_Generalization_of} consider parametric minimum cut problems where the capacities of the arcs incident to one specific vertex~$u^* \neq s,t$ are non-decreasing functions of a parameter $\lambda$, whereas all other arc capacities are constant.
They prove that, for such instances, the minimum cuts satisfy the generalized nesting property, which implies that the optimal cost curve has at most $n-1$~breakpoints.

\medskip

\citeauthor{McCormick_1999_Fast_Algorithms}~\cite{McCormick_1999_Fast_Algorithms} extends the results of \citeauthor{Gallo_1989_A_Fast}~\cite{Gallo_1989_A_Fast} to more general multi-parametric minimum cut problems with piece-wise linear capacity functions and additional arcs not incident to~$s$ or~$t$ with parameter-dependent capacities.
By again applying an extension of the push-relabel algorithm, the produced minimum cuts satisfy the nesting property, so the optimal cost curve has at most $n-1$ breakpoints.

\medskip

There is extensive literature that builds on these three fundamental papers and generalizes their theory. Extensions of \citeauthor{Gallo_1989_A_Fast}~\cite{Gallo_1989_A_Fast} include, e.g., the $\mathcal{O}(mn\log(\nicefrac{n^2}{m}))$ maximum distance variant of push-relabel due to \citeauthor{Gusfield_1994_A_faster}~\cite{Gusfield_1994_A_faster}, \citeauthor{Hochbaum_2008_The_Pseudoflow}'s~\cite{Hochbaum_2008_The_Pseudoflow} $\mathcal{O}(mn\log(n))$ pseudoflow algorithm, and the divide-and-conquer algorithm from \citeauthor{Tarjan_2006_Balancing_Applied}~\cite{Tarjan_2006_Balancing_Applied}.

\citeauthor{Brumelle_2005_Ordered_optimal}~\cite{Brumelle_2005_Ordered_optimal} comprises both the cases investigated by \citeauthor{Gallo_1989_A_Fast}~\cite{Gallo_1989_A_Fast} and by \citeauthor{Arai_1993_Generalization_of}~\cite{Arai_1993_Generalization_of}, and partially the case studied by \citeauthor{McCormick_1999_Fast_Algorithms}~\cite{McCormick_1999_Fast_Algorithms}.
They present an algebraic sufficient condition for the existence of a selection of optimal solutions in a parametric optimization problem that are totally ordered.
However, the order of the critical regions that belong to the optimal solutions in this selection does not necessarily adhere to the same ordering.

\citeauthor{Scutella_2007_A_note}~\cite{Scutella_2007_A_note} extends the results of \citeauthor{Gallo_1989_A_Fast}~\cite{Gallo_1989_A_Fast} and \citeauthor{Arai_1993_Generalization_of}~\cite{Arai_1993_Generalization_of} to networks in which the parametrization of the arc capacities can involve both the source and the target and one additional vertex~$u^*$.
She shows that the minimum cuts for such instances satisfy a relaxed form of the generalized nesting property, which again implies that this special case of a single-parametric minimum cut problem has at most $n-1$~breakpoints.
Finding these breakpoints needs asymptotically as much time as a single maximum flow computation and is based on an improvement of the maximum flow algorithm by \citeauthor{Goldberg_1988_A_new}~\cite{Goldberg_1988_A_new}.

\citeauthor{Granot_2012_Structural_and}~\cite{Granot_2012_Structural_and} identify three larger classes of parametric minimum cut problems satisfying the nesting property.
Again, all minimum cuts can be computed in the same asymptotic time as a single minimum cut by using a specialized flow update strategy to move from one value of the parameter to the next.

\medskip

If we turn the parametric problem into the related parametric search problem, there is further literature (see~\cite{Granot_2012_Structural_and}).
For example, the \emph{discrete Newton algorithm}~\cite{McCormick_1994_Computing_maximum,Radzik_1992_Minimizing_capacity,Radzik_1992_Newtons_method,Radzik_1993_Parametric_Flows,Radzik_1998_Fractional_Combinatorial} efficiently solves the parametric search problem and can be extended to other parametric problems.

\medskip

The literature covered so far only covers single-parametric minimum cut problems.
\citeauthor{Aissi_2015_Strongly_polynomial}~\cite{Aissi_2015_Strongly_polynomial} show that, in the $p$-parametric \emph{global} minimum cut problem, the number of cuts that become minimum at some point in the parameter space is in $\mathcal{O}(m^pn^2 \log^{p-1}(n))$.
For the special case of two parameters, they give a tighter bound of $\mathcal{O}(n^3\log(n))$ cuts that can be enumerated in $\tilde{\mathcal{O}}(m^2n^4)$ time.

\citeauthor{Karger_2016_Enumerating_parametric}~\cite{Karger_2016_Enumerating_parametric} requires stronger conditions than \citeauthor{Aissi_2015_Strongly_polynomial}~\cite{Aissi_2015_Strongly_polynomial} to apply his results.
The conditions are that every~$\lambda\in\Lambda$ and~$c_0(a),\ldots,c_p(a)$ for all~$a\in A$ are nonnegative.
Under these conditions, the results are stronger and can be generalized in different directions.
For instance, the $\mathcal{O}(n^{p+1})$ $p$-parametric \emph{global} minimum cuts can be enumerated in $\tilde{\mathcal{O}}(mn^{p+1})$ time.
He also gives a first generalization to single-parametric \emph{global} minimum cuts with polynomial cost functions.

\subsection{Parametric Minimum Cost Flow Problems}\label{subsec::MinCostFlow}

In the parametric \emph{minimum cost flow problem} (MCF), a graph $G=(V,A)$ with $|V|=n$ vertices and $|A|=m$ arcs together with a source vertex~$s \in V$, a target vertex~$t \in V$ with $s \neq t$, a demand~$D$, and a parameter set~$\Lambda \subseteq \mathbb{R}^p$ are given.
Additionally, for every arc~$a\in A$, a lower capacity~$l_a$, an upper capacity~$u_a$, and a parametric arc cost function~$c_a(\lambda)=c_0(a)+\sum_{i=1}^p\lambda\cdot c_i(a)$ are given.
Every feasible solution~$x\colon A\to\mathbb{R}$ represents an $s$-$t$-flow with flow value~$D$ that satisfies the capacity constraints $l_a\leq x(a) \leq u_a$ for all~$a\in A$.
The goal is to compute, for every $\lambda \in \Lambda$, a feasible $s$-$t$-flow $x_\lambda$ that minimizes $\sum_{a\in A} x_{\lambda}(a)\cdot c_a(\lambda)$.

\medskip

Based on a graph from \citeauthor{Zadeh_1973_A_bad}~\cite{Zadeh_1973_A_bad}, \citeauthor{Ruhe_1988_Complexity_results}~\cite{Ruhe_1988_Complexity_results} shows that even single-parametric MCF instances may have optimal solution sets of exponential size.

\subsubsection{Exact Methods}

For the MCF, every solution that is efficient in the multi-objective sense is also optimal for at least one parameter value in the parametric problem.
Hence, many important results from the multi-objective literature are directly applicable to the parametric setting.
\citeauthor{Hamacher_2007_Multiple_objective}~\cite{Hamacher_2007_Multiple_objective} provide a comprehensive survey of bi-objective MCF algorithms.
Many of these algorithms can be used for the single-parametric MCF as well.
We refer to some important algorithms described therein.
For detailed information on these, we refer directly to \citeauthor{Hamacher_2007_Multiple_objective}~\cite{Hamacher_2007_Multiple_objective}.

\medskip

Network simplex algorithms can be applied to the parametric MCF.
The idea is to pivot between basic feasible flows.
Such an algorithm is described, for example, by \citeauthor{pulat.huarng.ea1992EfficientSolutionsBicriteria}~\cite{pulat.huarng.ea1992EfficientSolutionsBicriteria}.
\citeauthor{lee.siminpulat1991BicriteriaNetworkFlow}~\cite{lee.siminpulat1991BicriteriaNetworkFlow} develop an algorithm that combines an out-of-kilter algorithm and a parametric linear programming procedure.
Using these two approaches, \citeauthor{sedeno-noda.gonzalez-martin2000BiobjectiveMinimumCost}~\cite{sedeno-noda.gonzalez-martin2000BiobjectiveMinimumCost} provide an algorithm that outperforms them both in a computational study.
Another algorithm by \citeauthor{sedeno-noda.gonzalez-martin2003AlternativeMethodSolve}~\cite{sedeno-noda.gonzalez-martin2003AlternativeMethodSolve} specializes the Eisner-Severance method for the parametric MCF, but is outperformed by their earlier algorithm from~\cite{sedeno-noda.gonzalez-martin2000BiobjectiveMinimumCost}.
A primal-dual network simplex method is developed by \citeauthor{eusebio.figueira.ea2009PrimalDualSimplex}~\cite{eusebio.figueira.ea2009PrimalDualSimplex}.
In a computational study, it outperforms a parametric simplex algorithm but is not compared to any of the aforementioned algorithms.
The network simplex algorithm from \citeauthor{raith.sedeno-noda2017FindingExtremeSupported}~\cite{raith.sedeno-noda2017FindingExtremeSupported} uses a labeling strategy to avoid computing redundant minimum ratios.
This leads to an output-polynomial running time in $\mathcal{O}(T_{\mathrm{MCF}}+ Bn(m+n\log n))$. In a computational study, the algorithm outperforms the Eisner-Severance method and a parametric simplex algorithm by \citeauthor{Raith_2009_A_two-phase}~\cite{Raith_2009_A_two-phase}.

\medskip

In the parametric literature, recent interest is mostly focused on the single-parametric MCF with parametric lower bound constraints.
Here, for every arc, the cost~$c_a(\lambda)$ is a constant function of~$\lambda$, but the lower capacities are of the form~$l_a(\lambda)=l_{0}(a)-\lambda\cdot l_{1}(a)$.

\medskip

Taking inspiration from generative linguistics and the behavior of neural networks, \citeauthor{Sangeorzan_2010_Partitioning_preflow-pull}~\cite{Sangeorzan_2010_Partitioning_preflow-pull} present an algorithm for the MCF with single-parametric lower bounds.
The algorithm recursively partitions the flows via a \emph{preflow-pull} strategy.
Its running time is in $\mathcal{O}(n^2m^{1/2}+Bn)$ and can be sped up by using parallelization strategies.

\medskip

In a series of papers, \citeauthor{Ciurea_2010_Balancing_algorithm}~\cite{Ciurea_2010_Balancing_algorithm, Ciurea_2012_A_sequential,Ciurea_2013_Shortest_conditional,Parpalea_2016_Minimum_Parametric,Parpalea_2019_Parametric_flows} develop several algorithms for the MCF with single-parametric lower bounds.
An overview of the most important algorithmic results is provided in~\cite{Parpalea_2019_Parametric_flows}.
The first algorithm, given by \citeauthor{Ciurea_2010_Balancing_algorithm}~\cite{Ciurea_2010_Balancing_algorithm}, covers only the special case of monotone parametric bipartite networks and is not always finite.
An algorithm for the general case is developed in~\cite{Ciurea_2012_A_sequential}.
It iteratively computes the breakpoints and a minimum cost flow for each critical region in ascending order of the parameter values.
Each iteration consists of modifying the network and then solving a maximum cost flow problem on it.
Based on the resulting flow, the next breakpoint and minimum cost flow can be computed.
It runs in time $\mathcal{O}(B \, T_{\mathrm{MaxFlow}})$, where $T_\mathrm{MaxFlow}$ denotes the time to solve a single maximum flow problem.

In~\cite{Ciurea_2013_Shortest_conditional}, \citeauthor{Ciurea_2013_Shortest_conditional} describe an algorithm that uses so-called \emph{shortest conditional decreasing paths} on a residual network.
The basic idea is to partition the parameter set and to repeatedly investigate each critical region and update it.
Each investigation consist of finding shortest conditional decreasing paths to determine if an update is needed.
The running time of this algorithm is in $\mathcal{O}(Bmn(B+mn))$.

A similar partitioning approach is taken in~\cite{Parpalea_2016_Minimum_Parametric}.
However, in the residual networks constructed there, it is sufficient to find any directed path since such a path must then necessarily be a conditional decreasing path as well.
This leads to a better running time in $\mathcal{O}(Bn^2m)$.

\subsection{Parametric Matroid and Spanning Tree Problems}\label{subsec::MinSpannTree}

Before we turn to the parametric matroid problem and some of its special cases, we define matroids as follows:

\begin{definition}[Independent Sets and Matroids]\label{Def:Matroid}
    A \emph{matroid} consists of a finite ground set~$S$ and a family~$\mathcal{I}$ of subsets of~$S$, called \emph{independent sets}, with the following properties:
    \begin{itemize}
        \item [I1)]  $\emptyset \in \mathcal{I}$
        \item [I2)] $\forall \, B\subseteq A\subseteq S:$ $A\in \mathcal {I}$ $\Longrightarrow \, B\in \mathcal {I}$
        \item [I3)] $\forall \, A,B \in \mathcal{I}:$ $|A|>|B|$ $ \Longrightarrow \exists \, x\in A \setminus B:$ $B\cup \{x\} \in \mathcal {I}$
    \end{itemize}

An inclusion-wise maximal independent set of a matroid~$(S,\mathcal{I})$ is called a \emph{basis}, and the \emph{rank} of a matroid is the cardinality of a basis.
\end{definition}

In the \emph{parametric matroid problem}, a matroid~$(S,\mathcal{I})$ with $|S|=n$, a parameter set $\Lambda \subseteq \mathbb{R}^p$, and parametric costs for the elements of~$S$ are given.
The feasible solutions correspond to the bases of the matroid, and the objective function value~$F_\lambda(x)$ of a basis~$x$ is given as the total cost of the elements of~$x$.
The goal is to obtain, for every~$\lambda\in\Lambda$, a matroid basis~$x_{\lambda}$ of minimum total cost.

\medskip

A special case of the parametric matroid problem is the parametric \emph{minimum spanning tree problem} (MST) (also called the \emph{parametric graphic matroid problem}).
Here, a connected graph $G=(V,E)$ with $|V|=n$ vertices, $|E|=m$ edges, a parameter set $\Lambda \subseteq \mathbb{R}^p$, and parametric edge costs are given. The ground set is the set~$E$ of edges of~$G$, and matroid bases correspond to spanning trees of~$G$. Thus, the goal is to find, for every~$\lambda\in\Lambda$, a spanning tree of minimum total edge cost.

\medskip

Another special case is the \emph{parametric uniform matroid problem}, where a ground set~$U_r$ of cardinality~$n$ is given, and matroid bases correspond to subsets containing exactly~$r$ elements. The goal is to find, for every~$\lambda\in\Lambda$, a matroid basis~$x_{\lambda}$ of minimum total cost. 

\subsubsection{Cardinality Bounds}
Parametric matroid and minimum spanning tree problems are well-researched topics in the parametric literature.
In contrast to many other problems, they always have an optimal solution set of polynomial cardinality in the single-parametric case.
The most important results for matroids are summarized in \Cref{Table:MatroidBounds}.

\medskip

\begin{table}[h]
    \footnotesize
    \centering
    \begin{tabular}{lp{3cm}p{3.5cm}l}
        \toprule
        \multicolumn{1}{l}{\textbf{Problem} ($p=1$)} & \multicolumn{1}{l}{\textbf{Setting \& Notation}} & \multicolumn{1}{l}{\textbf{Cardinality Bound}} & \multicolumn{1}{l}{\textbf{Reference}} \\
        \midrule[2pt]
        General Matroid & $n$ elements, rank $k$ & $\Omega(nk^{1/3})$ & \cite{Eppstein_1998_Geometric_Lower} \\
        \addlinespace
        & & $\mathcal{O}(n\cdot \min(k,n-k)^{1/2})$ & \cite{Gusfield_1980_Bounds_for} \\
        \midrule
        Graphic Matroid & $|V|=n$, $|E|=m$ & $\Omega(m\alpha(n))$ & \cite{Eppstein_1998_Geometric_Lower} \\ \addlinespace
        & & $\Omega(m\log(n))$ & \cite{Eppstein_2023_A_Stronger} \\
        \addlinespace
        & & $\mathcal{O}(m\sqrt{n})$ & \cite{Gusfield_1980_Bounds_for} \\ \addlinespace
        & & $\mathcal{O}(mn^{1/3})$ & \cite{Dey_1998_Improved_Bounds} \\
        \midrule
        Uniform Matroid & $n$ elements, rank $k$ & $\Omega(n\log(k))$ & \cite{Erdös_1973_Dissection_Graphs} \\ \addlinespace
        & & $\mathcal{O}(nk^{1/2})$ & \cite{Lovasz_1971_On_the} \\ \addlinespace
        & & $\mathcal{O}(\nicefrac{nk^{1/2}}{\log^*(k)})$ & \cite{Pach_1992_An_upper} \\ \addlinespace
        & & $\mathcal{O}(nk^{1/3})$ & \cite{Dey_1998_Improved_Bounds} \\
        \bottomrule
    \end{tabular}
    \caption{Upper and lower cardinality bounds for single-parametric matroid and minimum spanning tree problems. The function $\alpha$ denotes the inverse Ackermann function.}\label{Table:MatroidBounds}
\end{table} 

\medskip

For the parametric MST, the results of \citeauthor{Seipp_2013_On_adjacency}~\cite{Seipp_2013_On_adjacency} for the multi-objective MST can be transferred to the parametric setting.
They imply that, for every fixed number~$p\geq1$ of parameters, an optimal solution set of the parametric MST has a cardinality in $\mathcal{O}(m^{2(p+1)})$.

\medskip

If property I3) is removed from Definition \ref{Def:Matroid}, the resulting problem is called \emph{parametric independence systems problem}.
Because intractable problems such as the parametric maximum matching problem and the parametric assignment problem are special cases of it, optimal solution sets can have a cardinality exponential in $n$~\cite{Helfrich_2022_An_approximation}.

\subsubsection{Exact Methods for the Spanning Tree Problem}

\citeauthor{Fernandez-Baca_1996_Using_sparsification}~\cite{Fernandez-Baca_1996_Using_sparsification} provide an exact algorithm that solves the single-parametric MST in time $\mathcal{O}(\min\{nm\log(n),T_{\text{MST}}(2n,n)\cdot B\})$, where $T_{\text{MST}}(m,n)$ denotes the time for computing a minimum spanning tree in a graph with $m$ edges and $n$ vertices.
Their main idea is sparsification in order to pass through a reduced number of potential breakpoints.
This gives a faster asymptotic running time than the Eisner-Severance method~\cite{Fernandez-Baca_1996_Using_sparsification}.

One motivation for improving the running times of algorithms for the parametric MST is the wide range of applications.
Note that the construction procedure of \citeauthor{Fernandez-Baca_1996_Using_sparsification}~\cite{Fernandez-Baca_1996_Using_sparsification} can, among others~\cite{Hassin_1989_Maximizing_Classes}, also be used for stochastic spanning trees~\cite{Ishii_1981_Stochastic_spanning}.
Here, edge costs are not constant but random variables and the objective is to find an optimal spanning tree satisfying a certain probability constraint.
After transforming the stochastic problem into the deterministic equivalent problem and including a positive parameter in the objective function, their parametric algorithm finds an optimal solution in $\mathcal{O}(n^6)$ time.

\medskip

\citeauthor{Agarwal_1998_Parametric_and}~\cite{Agarwal_1998_Parametric_and} develop deterministic and randomized algorithms for the single-parametric MST and the kinetic MST on general graphs as well as improved algorithms for special families such as minor-closed and planar graphs.
In the \emph{(structurally) kinetic minimum spanning tree problem}, edges can be inserted or removed by parameter changes.
A multitude of well-known techniques are used to develop the algorithms therein, for example graph sparsification, convex hulls, and parametric search.
For general graphs, algorithms with running time in $\mathcal{O}(n^{2/3}\log^{4/3}(n))$ (or randomized in $\mathcal{O}(n^{2/3}\log(n))$) are presented.
If the graph is planar or from other minor-closed families of graphs, the running time bound for the parametric MST reduces to $\mathcal{O}(n^{1/4}\log^{3/2}(n))$ (or randomized to $\mathcal{O}(n^{1/4}\log(n))$), and for the kinetic MST to $\mathcal{O}(n^{1/2}\log^{3/2}(n))$ (or randomized to $\mathcal{O}(n^{1/2}\log(n))$).

\subsubsection{Algorithms for Different Settings}

In the multi-objective literature, an algorithm for the parametric MST is described by \citeauthor{Seipp_2013_On_adjacency}~\cite{Seipp_2013_On_adjacency}.
His approach works by computing the minimization diagram (see \Cref{subsec::decomposition}) in two steps:
First, it uses some special properties of the parametric MST to find an arrangement of hyperplanes that subdivides~$\Lambda$ into polyhedra.
Then, it iteratively checks for adjacent polyhedra whether they are part of the same critical region and merges such polyhedra until only the critical regions of a  minimization diagram remain.
No exact running time analysis is given, but a polynomial running time bound is shown.

\medskip

There is also a lot of literature on parametric search or finding the shortest bottleneck edge for the parametric MST.
However, this is not the focus of this survey and related literature can be found, for example, in~\cite{Fernandez-Baca_1997_Linear-time_algorithms,Fernandez-Baca_2000_Multi-parameter_Minimum,Katoh_2001_Notes_on,Chan_2005_Finding_the}.

Furthermore, there are other parametric problems related to the main problems discussed in this section. For example, \citeauthor{Hausbrandt__2024_Parametric_matroid}~\cite{Hausbrandt__2024_Parametric_matroid} develops upper and lower bounds as well as algorithms for the parametric matroid interdiction problem.
Due to the limited literature available, these problems are also not in our focus.

\subsection{Parametric Knapsack Problems}\label{subsec::Knapsack}

The \emph{parametric knapsack problem} (sometimes also called the \emph{parametric $0$-$1$-knapsack problem}) is defined by a item weight vector $w = (w_1,\dots,w_n)^\mathsf{T} \in \mathbb{N}^n$, a budget $K\in\mathbb{N}$, a parameter set $\Lambda \subseteq \mathbb{R}^p$, and a parametric item profit vector $c(\lambda) = (c_1(\lambda),\dots,c_n(\lambda))^\mathsf{T} \in \mathbb{R}^n$ with $c_i(\lambda)=c_0^i + \sum^p_{j=1}\lambda_j\cdot c_j^i$ for $i \in [n]$ and $\lambda \in \Lambda$.
A feasible solution $x\in\{0,1\}^n$ (a \emph{packing}) represents a subset of selected items such that $w^\mathsf{T}x\leq K$.
The goal is to find, for every~$\lambda\in\Lambda$, a feasible packing~$x \in\{0,1\}^n$ maximizing the total profit $c(\lambda)^\mathsf{T}x$ of the selected items.

\subsubsection{Intractability Results and Exact Methods}
\citeauthor{Carstensen_1983_Complexity_of}~\cite{Carstensen_1983_Complexity_of} is the first to address the parametric knapsack problem and shows that the cardinality of an optimal solution set for a single-parametric knapsack problem can be exponential in the number of items.
This result even holds under the additional assumption that~$\Lambda$ is a compact interval in~$\mathbb{R}_{>0}$ and that $c_0^i,c_1^i>0$ for all~$i \in [n]$.
 \citeauthor{Carstensen_1983_Complexity_of}~\cite{Carstensen_1983_Complexity_of} also shows, however, that the cardinality of an optimal solution set of any single-parametric \emph{binary} integer program is linear in the number of variables when the cost function coefficients~$c_0^{i}$ and~$c_1^{i}$ are integers in $[-M,M]$ for some constant~$M \in \mathbb{N}$ and all~$i \in [n]$.
In particular, this implies that the cardinality of an optimal solution set of the single-parametric knapsack problem is linear in the number of items under this condition.

\medskip

An algorithm for the single-parametric knapsack problem is developed by \citeauthor{Eben-Chaime_1996_Parametric_Solution}~\cite{Eben-Chaime_1996_Parametric_Solution}. This algorithm combines the Eisner-Severance method with a classic dynamic programming approach for knapsack problems.
Its running time is in~$\mathcal{O}(KnB)$.

\medskip

A problem closely related to the parametric knapsack problem is the \emph{inverse-parametric knapsack problem}, where
the goal is to find the smallest parameter for which the optimal objective value of a parametric knapsack problem is equal to a specific input value.
For this problem, \citeauthor{Burkard_1995_The_inverse-parametric}~\cite{Burkard_1995_The_inverse-parametric} provide a pseudo-polynomial algorithm based on Megiddo's method with a running time in $\mathcal{O} \big(\nicefrac{n^2 K \log(K)}{\log(n\log(K))}\big)$. \medskip

\subsubsection{Parametric Approximation Methods}

Several approximation algorithms exist for the single-parametric knapsack problem.
\citeauthor{Giudici_2017_Approximation_schemes}~\cite{Giudici_2017_Approximation_schemes} provide a PTAS for the single-parametric knapsack problem.
The authors also describe how approximation algorithms of bi-objective optimization problems can be used for approximating single-parametric optimization problems.
They use this result and the algorithm from \citeauthor{Erlebach_2002_Approximating_multiobjective}~\cite{Erlebach_2002_Approximating_multiobjective} to derive an FPTAS for the case that~$\Lambda=\mathbb{R}_{\geq0}$ and all profit coefficients in the objective function are nonnegative.

\citeauthor{Holzhauser_2017_An_FPTAS}~\cite{Holzhauser_2017_An_FPTAS} describe an FPTAS that works on arbitrary parameter intervals and profit coefficients in the objective function.
Furthermore, in contrast to the FPTAS from \citeauthor{Giudici_2017_Approximation_schemes}~\cite{Giudici_2017_Approximation_schemes}, the running time is strongly polynomial.
For any~$\varepsilon>0$, the algorithm finds a $(1-\varepsilon)$-approximation set with cardinality in $\mathcal{O}\big(\nicefrac{n^2}{\varepsilon}\big)$ and has a running time in $\mathcal{O}\big(\nicefrac{n^2}{\varepsilon}\,\cdot\, T_\Pi(\varepsilon)\big)$, where $T_\Pi(\varepsilon)$ denotes the running time of an FPTAS for the conventional knapsack problem with an approximation factor $1+\varepsilon$.

\medskip

\citeauthor{Halman_2018_An_FPTAS}~\cite{Halman_2018_An_FPTAS} investigate a variant of the single-parametric knapsack problem in which the parameter dependence is in the weight vector instead of the objective function.
They provide two FPTASs, one of which has a running time in $\mathcal{O}\big( \nicefrac{n^3}{\varepsilon^2} \cdot \min\{\log^2(P),n^2\} \cdot \min\{\log(M), \nicefrac{n \log(\nicefrac{n}{\varepsilon})}{\log(n \log(\nicefrac{n}{\varepsilon}))}\} \big)$, where $\varepsilon >0$ determines the approximation precision, $P$ is an upper bound on the optimal profit, and~$M$ is an upper bound on the reachable weight.

\subsection{Parametric Matching and Assignment Problems}\label{subsec::Matching}

In the \emph{parametric assignment problem}, a bipartite, undirected graph $G = (U \cup V , E)$ with $|U |=|V |=n$ vertices and $|E|=m$ edges, a parameter set $\Lambda \subseteq \mathbb{R}^p$, and parametric edge weights $w_{e}(\lambda)=w_0(e) + \sum_{i=1}^{p}\lambda_i \cdot w_i(e)$ for~$\lambda\in\Lambda$ and~$e \in E$ are given.
An assignment is a set of edges where every edge is between $U$ and $V$, and no vertex from $U$ or $V$ is an endpoint more than once.
The goal of the parametric assignment problem is to compute, for every $\lambda \in \Lambda$, a minimum weight assignment.

\medskip

The parametric assignment problem is intractable, even in the single-parametric case.
To prove this, \citeauthor{Bazgan_2022_An_approximation}~\cite{Bazgan_2022_An_approximation} sketch how to use the non-parametric reduction of the shortest path problem to the assignment problem from \citeauthor{Lawler_2001_Combinatorial_optimization}~\cite{Lawler_2001_Combinatorial_optimization}.
The first proof of intractability is given by \citeauthor{Carstensen_1983_The_Complexity}~\cite{Carstensen_1983_The_Complexity} (referring to the reduction in~\cite{Hoffmann_1963_A_note}).

However, for the special case where the parameter coefficients $w_1(e)$ are either $0$ or $1$ for all $e \in E$ the single-parametric assignment problem becomes tractable.
A polynomial-time algorithm is given by \citeauthor{Gassner_2010_A_fast}~\cite{Gassner_2010_A_fast}. It repeatedly uses the parametric shortest path algorithm of \citeauthor{Young_1991_Faster_parametric}~\cite{Young_1991_Faster_parametric} and its running time is in $\mathcal{O}(mn+n^2\log(n))$.
\citeauthor{Rosenmann_2022_Computing_the}~\cite{Rosenmann_2022_Computing_the} uses the parametric algorithm from \citeauthor{Gassner_2010_A_fast}~\cite{Gassner_2010_A_fast} for computing the sequence of all $k$-cardinality assignments. \medskip

The parametric assignment problem can be seen as a special case of the parametric maximum matching problem.
Hence, the (single-)parametric maximum matching problem is also intractable~\cite{Helfrich_2022_An_approximation}.
In the literature, the parametric maximum matching problem is rarely investigated.
\citeauthor{Weber_1981_Sensitivity_analysis}~\cite{Weber_1981_Sensitivity_analysis} conducts a sensitivity analysis for the case where one edge weight depends on a single parameter.
\citeauthor{Gusfield_1989_Parametric_stable}~\cite{Gusfield_1989_Parametric_stable} consider the single-parametric stable matching problem and develop an efficient algorithm by reducing it to the parametric minimum cut problem and using the method of \citeauthor{Gallo_1989_A_Fast}~\cite{Gallo_1989_A_Fast}.


\section{Conclusion}\label{subsec::Conclusion}

\begin{figure}[h]
  \centering
  \includegraphics[width=\linewidth]{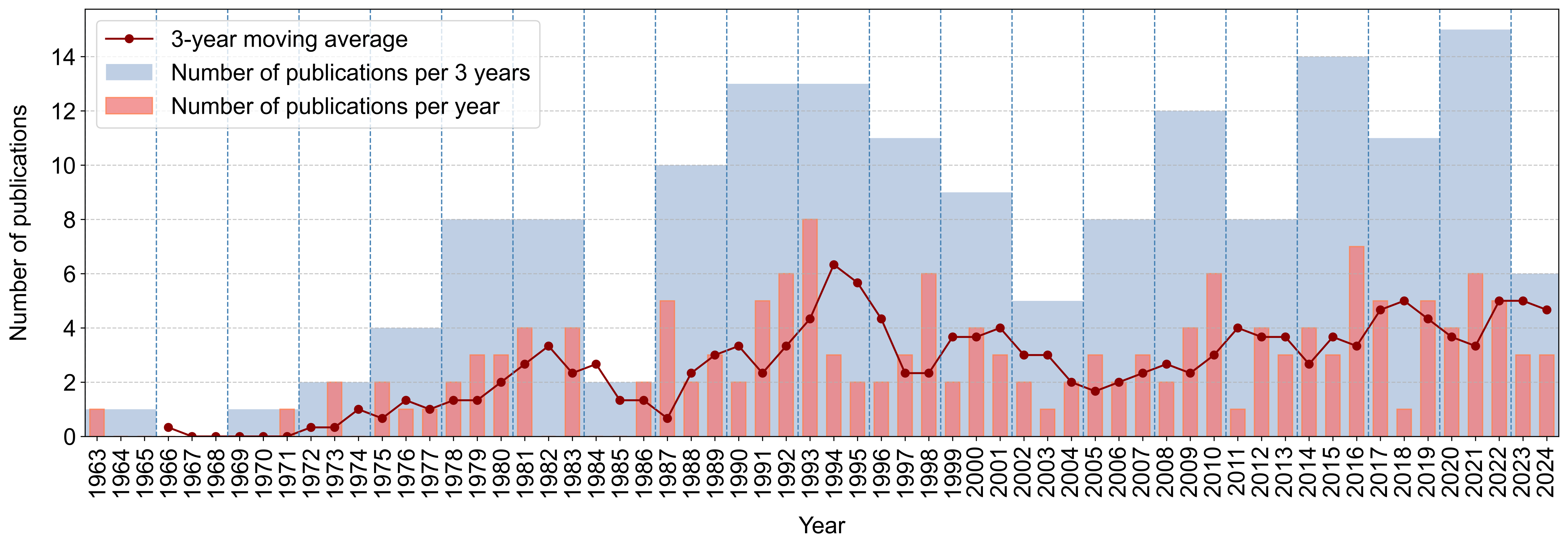}
  \caption{Number of publications covered in this survey over time (per year, per 3 years, and 3-year moving average) from 1963 to 2024. The moving
average for year $t$ considers the years~$t-1$, $t-2$, and~$t-3$.\label{fig::histogramm}}
\end{figure}

The literature on parametric optimization dates back over 60 years, and some contributions appear even earlier under different names such as sensitivity analysis~\cite{Gal_1997_Advances_in}.
\Cref{fig::histogramm} visualizes the distribution of the literature in this paper over the years.
Although it does not include papers that are solely concerned with practical applications, have somewhat different problem settings, or are simply surveyed in other literature reviews, \Cref{fig::histogramm} shows the ongoing interest in parametric optimization.

\medskip

In nearly all papers we surveyed, the parameter dependencies are linear and appear in the objective function.
Parameter-dependent constraints are mainly discussed for the generic MILP case.
But for specific problems, there are often structural properties that can be used to find cardinality bounds or develop specialized algorithms.
This can be seen in some of the algorithms we surveyed, for example, for minimum cost flow and maximum flow problems or for the knapsack problem.
For most other problems, however, such algorithms are rare.

\medskip

The main challenge in parametric optimization is intractability and finding strategies to cope with it.
Most parametric problems are intractable because of the cardinality of optimal solution sets.
This can be the case even if the underlying non-parametric optimization problem is solvable in polynomial time, which is the case, for example, for parametric shortest path problems.
In \Cref{fig::TableBig}, we summarize known bounds on the cardinality of optimal solution sets for all the specific problems considered in this survey.

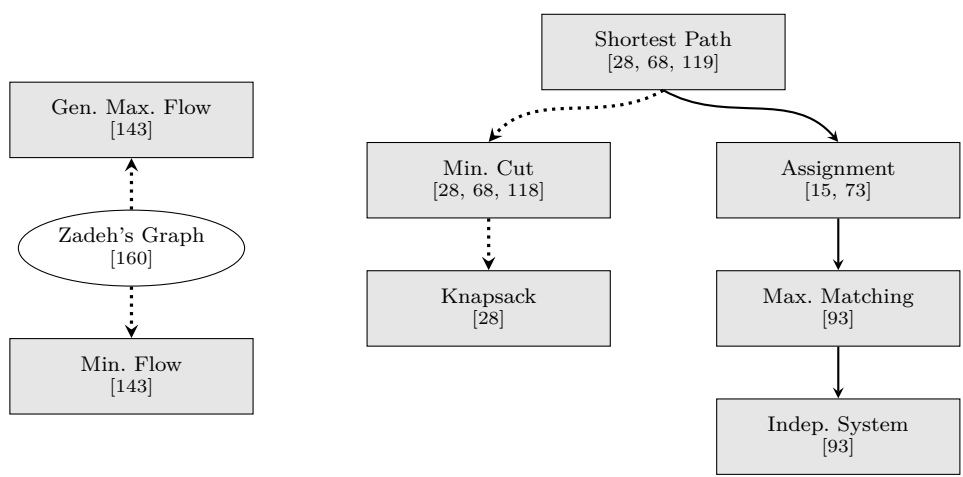
\begin{figure}[h]
  \centering
  \begin{tikzpicture}
    \def\xshift{2.3}
    \def\yshift{1.7}
    \tikzstyle{problem} = [rectangle, minimum width=3.2cm, minimum height=1cm, text centered, draw=black, fill=black!10, inner sep=4, align=center, font = {\footnotesize}]
    \tikzstyle{zadeh} = [ellipse, text centered, draw=black, fill=black!10, inner sep=2.5, align=center, font = {\footnotesize}];
    \tikzstyle{reduc} = [thick,->,>=stealth]
    \tikzstyle{construc} = [very thick,->,>=stealth, dotted]

    \node (shortpath) [problem] {{\small Shortest Path}\\\cite{Carstensen_1983_Complexity_of,Mulmuley_2001_A_Lower,Gajjar_2019_Parametric_Shortest}
    };
    \node (assign) [problem, xshift=\xshift cm, yshift=-\yshift cm] at (shortpath) {{\small Assignment}\\\cite{Gassner_2010_A_fast,Bazgan_2022_An_approximation}};
    \node (maxmatching) [problem, yshift=-\yshift cm] at (assign) {{\small Max.\@ Matching}\\\cite{Helfrich_2022_An_approximation}};
    \node (indsys) [problem, yshift=-\yshift cm] at (maxmatching) {{\small Indep.\@ System}\\\cite{Helfrich_2022_An_approximation}};

    \node (mincut) [problem, xshift=-\xshift cm, yshift=-\yshift cm] at (shortpath) {{\small Min.\@ Cut}\\\cite{Carstensen_1983_Complexity_of,Mulmuley_1999_Lower_bounds,Gajjar_2019_Parametric_Shortest}};
    \node (knapsack) [problem, yshift=-\yshift cm] at (mincut) {{\small Knapsack}\\\cite{Carstensen_1983_Complexity_of}};

    \node (zadehgraph) [zadeh, xshift=-7 cm, yshift=-2.6 cm,fill=white] at (shortpath) {{\small Zadeh's Graph}\\\cite{Zadeh_1973_A_bad}};
    \node (genmaxflow) [problem, xshift=0 cm, yshift=\yshift cm] at (zadehgraph) {{\small Gen.\@ Max.\@ Flow}\\\cite{Ruhe_1988_Complexity_results}};
    \node (minflow) [problem, xshift=0 cm, yshift=-\yshift cm] at (zadehgraph) {{\small Min.\@ Flow}\\\cite{Ruhe_1988_Complexity_results}};

    \draw[reduc, out=-30, in=130] (shortpath.south) to (assign.north);
    \draw[reduc] (assign.south) to (maxmatching.north);
    \draw[reduc] (maxmatching.south) to (indsys.north);

    \draw[construc, out=-150, in=50] (shortpath.south) to (mincut.north);
    \draw[construc] (mincut.south) to (knapsack.north);

    \draw[construc, in=-90, out=90] (zadehgraph.90) to (genmaxflow.south);
    \draw[construc, in=90, out=-90] (zadehgraph.-90) to (minflow.north);
\end{tikzpicture}
  \caption{Relations between intractability results from \Cref{subsec::specificProblems}. A solid arc indicates a proof via reduction, a dotted arc indicates that 
  constructions leading to optimal solution sets with superpolynomial cardinality are transferred into a different problem setting. \label{fig::IntractDiagram}}
\end{figure}

For proving intractability, many papers use similar strategies.
This is visualized in \Cref{fig::IntractDiagram}:
Either intractable parametric problems are reduced to the problem for which intractability is to be shown, or constructions that lead to optimal solution sets with superpolynomial cardinality are transferred into a different problem setting.

Still, some problems admit an optimal solution set of polynomial size.
This has been shown mainly for matroid-based problems.
Restricting the parameter coefficients to specific integers also often makes a problem tractable.
Problems where this can be observed are shortest path problems, the assignment problem, and binary integer programs with bounded cost function coefficients.

When dealing with intractable problems, a way to still achieve polynomial running time is to find an approximation set instead of an optimal solution set.
However, approximation algorithms are rarely found in the literature.
We discussed some general approximation algorithms in \Cref{subsec::generalAPX}.
To the best of our knowledge, the only problem-specific approximation algorithms can be found for the knapsack problem.
This demonstrates a gap in the literature -- particularly compared to multi-objective optimization, where problem-specific approximation algorithms are available for a large variety of problems~\cite{herzel.ruzika.ea2021ApproximationMethodsMultiobjective}.

With respect to methods for computing an optimal solution set, algorithms for specific problems are often adaptions of general strategies.
For example, the algorithms of \citeauthor{sedeno-noda.gonzalez-martin2003AlternativeMethodSolve}~\cite{sedeno-noda.gonzalez-martin2003AlternativeMethodSolve} for the parametric minimum cost flow problem, of \citeauthor{Fernandez-Baca_1996_Using_sparsification}~\cite{Fernandez-Baca_1996_Using_sparsification} for the parametric spanning tree problem, and of \citeauthor{Eben-Chaime_1996_Parametric_Solution}~\cite{Eben-Chaime_1996_Parametric_Solution} for the knapsack problem are tailored variants of the Eisner-Severance method, while 
Megiddo's method is used by \citeauthor{Burkard_1995_The_inverse-parametric}~\cite{Burkard_1995_The_inverse-parametric} for the parametric knapsack problem and by \citeauthor{Agarwal_1998_Parametric_and}~\cite{Agarwal_1998_Parametric_and} for the parametric spanning tree problem.
These adaption strategies further highlight the value of finding problem-agnostic methods.
Here, research areas related to parametric optimization may offer new opportunities.
For example, the immediate connection between parametric optimization and multi-objective optimization is used numerous times in this survey to transfer (in-)tractability results and algorithms, but similar connections might also exist to other research areas.

\begin{sidewaystable}[h]
    \footnotesize
    \resizebox{\linewidth}{!}{
    \begin{tabular}{lp{0.7cm}p{5cm}p{4cm}p{1cm}@{\hskip 1cm}p{6cm}}
        \toprule
        \multicolumn{2}{c}{\textbf{Parametric Problem}} & \multicolumn{3}{c}{\textbf{Cardinality Bounds}} & \multicolumn{1}{l}{\textbf{Algorithmic Results}}\\
        \cmidrule(lr){1-2}
        \cmidrule(lr){3-5}
        \cmidrule(lr){6-6}
        \multicolumn{1}{l}{\textbf{Problem}} & \multicolumn{1}{l}{\textbf{$p$}} & \multicolumn{1}{l}{\textbf{Setting \& Notation}} & \multicolumn{1}{l}{\textbf{Cardinality Bound}} & \multicolumn{1}{l}{\textbf{Reference}} & \multicolumn{1}{l}{\textbf{Running Time}} \\
        \midrule[2pt]
        $0$-$1$-Knapsack & $1$ & $n$ items & $2^n-1$ (lower bound) & \cite{Carstensen_1983_Complexity_of} & \,\,\,\,\multirow{3}{*}{\parbox{6cm}{$\mathcal{O}(KnB)$ for $K$~knapsack budget, $n$~items, $B$~breakpoints \cite{Eben-Chaime_1996_Parametric_Solution}}} \\ \addlinespace
        & $1$ & $n$ items, $f_i$ linear with integer weights bounded by $-C_i$ and $C_i$ & $\min\{(2C_1+1)n,(4C_0+2)n\}$ (upper bound) & \cite{Carstensen_1983_Complexity_of} & \\ \midrule
        Assignment & $1$ & $|V|+|U|=2n$ vertices & $2^{\Omega(\log^2(n))}$ & \cite{Bazgan_2022_An_approximation} & \,\,\,\,\multirow{2}{*}{\parbox{6cm}{$\mathcal{O}(mn+n^2\log n)$ for special case $F_{\lambda}(x)=\sum_{e\in E}(a_e+\lambda b_e)\cdot x_e$ with $b_e \in \{-1,\,0\} \, \forall e \in E$\cite{Gassner_2010_A_fast}}}\\ \addlinespace
        & $1$ & $|V|+|U|=2n$ vertices & $2^{\mathcal{O}(n)}$ & \cite{Hrubes_2023_Shadows_of} & \\ \addlinespace
        \midrule
        Min. Cut & $1$ & $n+2$ vertices & $2^n-1$ (lower bound) & \cite{Carstensen_1983_Complexity_of} & \,\,\,\,\multirow{4}{*}{\parbox{6cm}{$\Tilde{\mathcal{O}}(mn^{1+p})$ for special case of finding \textit{global} min. cuts in $p$ dimensions with nonneg. weights and parameters \cite{Karger_2016_Enumerating_parametric}}} \\ \addlinespace
        & $1$ & planar $G$, polynomial number vertices in~$n$ & $2^{\Omega(\log^2(n))}$ & \cite{Gajjar_2019_Parametric_Shortest} & \\ \addlinespace
        & $2$ & $n+2$ vertices, fulfilling SSM & $2^n-1$ (lower bound) & \cite{Allman_2022_Complexity_of} &
        \,\,\, \multirow{3}{*}{\parbox{6cm}{$\mathcal{O}(1)$ max. flow computations for special case of single-parametric problem fulfilling SSM \cite{Gallo_1989_A_Fast}}} \\ \addlinespace
        & $1$ & $n+2$ vertices, fulfilling SSM & $n-1$ (upper bound) & \cite{Eisner_1976_Mathematical_Techniques,Stone_1978_Critical_Load} & \\
        \midrule
        Shortest Path & $1$ & layered \& planar $G$, $n$ vertices & $2^{\Omega(\log^2(n))}$ & \cite{Gajjar_2019_Parametric_Shortest} & \,\,\,\,\multirow{2}{*}{\parbox{6cm}{$\mathcal{O}(mn+n^2\log n)$ for single-source case, $F_{\lambda}(x)=\sum_{e\in E}(a_e+\lambda b_e)\cdot x_e$, $b_e \in \{-1,$ $0\} \, \forall e \in E$ \cite{Young_1991_Faster_parametric}}} \\ \addlinespace
        & $1$ & $n$ vertices & $n^{\log{n}+\mathcal{O}(1)}$ & \cite{Carstensen_1983_The_Complexity} & \\ \addlinespace
        \midrule
        General Matroid & $1$ & $n$ elements, rank $k$ & $\mathcal{O}(n\cdot \min(k,n-k)^{1/2})$ & \cite{Gusfield_1980_Bounds_for} & \\ \addlinespace
        & $1$ & $n$ elements, rank $k$ & $\Omega(nk^{1/3})$ & \cite{Eppstein_1998_Geometric_Lower} & \\
        \midrule
        Spanning Tree & $1$ & $|V|=n$, $|E|=m$ & $\mathcal{O}(m\sqrt{n})$ & \cite{Gusfield_1980_Bounds_for} & \,\,\,\,\multirow{6}{*}{\parbox{6cm}{$\mathcal{O}(n^{1/4}\log^{3/2}(n)B)$ for planar graphs, $B$~breakpoints \cite{Agarwal_1998_Parametric_and}}}\\ \addlinespace
        & $1$ & $|V|=n$, $|E|=m$ & $\mathcal{O}(mn^{1/3})$ & \cite{Dey_1998_Improved_Bounds} & \\ \addlinespace
        & $1$ & $|V|=n$, $|E|=m$, inverse Ackermann function $\alpha$ & $\Omega(m\alpha(n))$ & \cite{Eppstein_1998_Geometric_Lower} & \\ \addlinespace
        & $1$ & $|V|=n$, $|E|=m$ & $\Omega(m\log n)$ & \cite{Eppstein_2023_A_Stronger} & \\ \addlinespace
        & $p-1$ & $|E|=m$ & $\mathcal{O}(m^{2p})$ & \cite{Seipp_2013_On_adjacency} & \\
        \midrule
        Uniform Matroid & $1$ & $n$ elements, rank $k$ & $\mathcal{O}(nk^{1/2})$ & \cite{Lovasz_1971_On_the} & \\ \addlinespace
        & $1$ & $n$ elements, rank $k$ & $\mathcal{O}(\nicefrac{nk^{1/2}}{\log^*(k)})$ & \cite{Pach_1992_An_upper} & \\ \addlinespace
        & $1$ & $n$ elements, rank $k$ & $\Omega(n\log(k))$ & \cite{Erdös_1973_Dissection_Graphs} & \\ \addlinespace
        & $1$ & $n$ elements, rank $k$ & $\mathcal{O}(nk^{1/3})$ & \cite{Dey_1998_Improved_Bounds} & \\
        \midrule
        Min.\ Cost Flow & $1$ & $n+2$ vertices & $2^n-1$ (lower bound) & \cite{Ruhe_1988_Complexity_results} & \,\,\,\,\multirow{2}{*}{\parbox{6cm}{$\mathcal{O}(n^2m^{1/2}+Bn)$ for $B$~breakpoints \cite{Sangeorzan_2010_Partitioning_preflow-pull}}} \\ \addlinespace[0.8cm]
        & & & & & \,\,\, \multirow{2}{*}{\parbox{6cm}{$\mathcal{O}(nm\log(n^2/m)\log(nC))$ for $C$ max. arc cost~\cite{Lin_2015_On_the}}} \\ \addlinespace[0.5cm]
        \midrule
        LP & $1$ & $x \in \mathbb{R}^{n_1}$, $y \in \mathbb{Z}^{n_2}$, $n_1$ even, ${m_1 = \nicefrac{n_1}{2}}$ for $m_1+m_2$ dim. of RHS & $2^{n_1/2}$ (lower bound) & \cite{Murty_1980_Computational_Complexity} & \\ \bottomrule
    \end{tabular} }
    \caption{Comparison of best upper/lower cardinality bounds and algorithmic results for specific parametric optimization problems. \label{fig::TableBig}}
\end{sidewaystable}

\FloatBarrier


\section*{Declarations}

\subsection*{Funding}
This research was funded by the Deutsche Forschungsgemeinschaft (DFG, German Research Foundation) -- Project number 508981269.

\subsection*{Author contributions}
\textbf{Levin Nemesch:} Conceptualization, Methodology, Investigation, Literature Review, Data Curation, Writing - Original Draft, Writing - Review and Editing, Visualization.
\textbf{Stefan Ruzika:} Conceptualization, Methodology, Writing - Review and Editing, Supervision, Project Administration, Funding acquisition.
\textbf{Clemens Thielen:} Conceptualization, Methodology, Writing - Review and Editing, Supervision, Project Administration, Funding acquisition.
\textbf{Alina Wittmann:} Conceptualization, Methodology, Investigation, Literature Review, Data Curation, Writing - Original Draft, Writing - Review and Editing, Visualization.


\printbibliography[heading=bibintoc]{}
\end{document}